\documentclass[11pt]{article}
\usepackage{amsmath,amsfonts,amssymb,amsthm,enumerate,graphicx}
\usepackage{tikz}
\usepackage{subfigure}
\usepackage{amsmath}
\usepackage{amssymb}
\usepackage[english]{babel}
\usepackage{graphicx}
\usepackage{gastex}
\usepackage{longtable}
\usepackage{lscape}
\usepackage{multicol}
\usepackage{latexsym}
\usepackage{verbatim}
\usepackage{multirow}
\usepackage{longtable}
\usepackage{authblk}
\usepackage{braket}
\usetikzlibrary{arrows}

\usepackage{amssymb,latexsym,amsmath,epsfig,amsthm,mathrsfs}
\usepackage{amsfonts}
\usepackage{amscd}
\usepackage{dsfont}
\usepackage{bbm}
\usepackage{rotating}
\usepackage{newlfont}
\usepackage{enumitem}
\usepackage{tabularx}
\usepackage{longtable}
\usepackage{tabularx,ragged2e,booktabs,caption}
\usepackage{multicol}
\usetikzlibrary{arrows}
\tikzstyle{vertex}=[ draw, inner sep=0pt, minimum size=0pt]

\usepackage{epstopdf}
\usepackage{breakurl}
\usepackage{hyperref}
\usepackage{float}
\usepackage{enumitem}
\usepackage{mathrsfs}
\usepackage[mathscr]{euscript}
\usepackage{mathtools}
\usepackage{tikz-cd}
\usepackage{array,multirow}
\usepackage[]{cite}
\usepackage{ragged2e}
\newlist{subquestion}{enumerate}{1}
\date{}
\newcommand{\circhash}{\mbox{ $\bigcirc\hspace{-3.5mm}\# $}}

\newtheorem{theorem}{{\bf Theorem}}[section]

\newtheorem{lemma}[theorem]{{\bf Lemma}}

\newtheorem{proposition}[theorem]{{\bf Proposition}}
\newtheorem{remark}[theorem]{{\bf Remark}}

\numberwithin{Subcase}{Case}

\numberwithin{Subsubcase}{Subcase}
\newtheorem{claim}{\textit{Claim}}

\newtheorem{cor}{{\sc Corollary}}[section]

\newcommand{\Echar}[1]{Euler characteristic #1}

\begin{document}
	
\title{Dihedral and cyclic symmetric maps on surfaces}

 	\author[1] {Marbarisha M. Kharkongor}
 	\author[2] {Debashis Bhowmik}
 	\author[1] {Dipendu Maity}
 	\affil[1]{Department of Sciences and Mathematics,
 		Indian Institute of Information Technology Guwahati, Bongora, Assam-781\,015, India.\linebreak
 		\{marbarisha.kharkongor, dipendu\}@iiitg.ac.in/\{marbarisha.kharkongor, dipendumaity\}@gmail.com.}
 	\affil[2]{Department of Mathematics, Indian Institute of Technology Patna, Patna 801\,106, India.
 		debashisiitg@gmail.com.}

\vspace{-15mm}

\date{\today}

\maketitle

\vspace{-10mm}

\begin{abstract}
If the face\mbox{-}cycles at all the vertices in a map are of the same type, then the map is said to be a semi-equivelar map. Automorphism (symmetry) of a map can be thought of as a permutation of the vertices which preserves the vertex\mbox{-}edge\mbox{-}face incidences in the embedding. The set of all symmetries forms the symmetry group. In this article, we discuss the maps' symmetric groups on higher genus surfaces. In particular, we show that there are at least $39$ types of the semi-equivelar maps on the surface with Euler char. $-2m, m \ge 2$ and the symmetry groups of the maps are isomorphic to the dihedral group or cyclic group. Further, we prove that these $39$ types of semi-equivelar maps are the only types on the surface with Euler char. $-2$.
Moreover, we know the complete list of semi-equivelar maps (up to isomorphism) for a few types. We extend this list to one more type and can classify others similarly. We skip this part in this article.   
\end{abstract}

\noindent {\small {\em MSC 2010\,:} 52C20, 52B70, 51M20, 57M60.

\noindent {\em Keywords:} Symmetric groups; Semi-equivelar maps; Polyhedral maps; Cyclic group; Dihedral group.}

\section{Introduction}
A map is a connected $2$-dimensional cell complex on a surface. Equivalently, it is a cellular embedding of a connected graph on a surface. In this article, a map will mean a polyhedral map on a surface, that is, non-empty intersection of any two faces is either a vertex or an edge. For a map $\mathcal{K}$, let $V(\mathcal{K})$ be the vertex set of $\mathcal{K}$ and $u\in V(\mathcal{K})$. The faces containing  $u $ form a cycle (called the {\em face-cycle} at  $u $)  $C_u $ in the dual graph  of  $\mathcal{K} $. That is,  $C_u $ is of the form  $(F_{1,1}\mbox{-}\cdots\mbox{-}F_{1,n_1})\mbox{-}\cdots\mbox{-}(F_{k,1}\mbox{-}\cdots \mbox{-}F_{k,n_k})\mbox{-}$ $F_{1,1} $, where  $F_{i,\ell} $ is a  $p_i $-gon for  $1\leq \ell \leq n_i $,  $1\leq i \leq k $,  $p_r\neq p_{r+1} $ for  $1\leq r\leq k-1 $ and  {$p_k\neq p_1 $}. A map  $\mathcal{K} $ is said to be {\em semi-equivelar} (see \cite{datta2017}) if  $C_u $ and  $C_v $ are of same type for all  $u, v \in V(\mathcal{K}) $. That is, there exist integers  $p_1, \dots, p_k\geq 3 $ and  $n_1, \dots, n_k\geq 1 $,  $p_i\neq p_{i+1} $ (addition in the suffix is modulo  $k $) such that  $C_u $ is of the form as above for all  $u\in V(\mathcal{K}) $. In such a case,  $\mathcal{K} $ is called a semi-equivelar map (SEM) of type  $[p_1^{n_1}, \dots, p_k^{n_k}] $ (or, a map of type  $[p_1^{n_1}, \dots, p_k^{n_k}] $). Let  $F_1\mbox{-}\cdots\mbox{-}F_m\mbox{-}F_1 $ be the face-cycle of a vertex  $u $ in a map $\mathcal{K}$. Then  $F_i \cap F_j $ is either  $u $ or an edge through  $u $. Thus the face  $F_i $ must be of the form  $u_{i+1}\mbox{-}u\mbox{-}u_i\mbox{-}P_i\mbox{-}u_{i+1} $, where  $P_i = \emptyset $ or a path  $\& $  $P_i \cap P_j = \emptyset $ for  $i \neq j $. Here addition in the suffix is modulo  $m $. So,  $u_{1}\mbox{-}P_1\mbox{-}u_2\mbox{-}\cdots\mbox{-}u_m\mbox{-}P_m\mbox{-}u_1 $ is a cycle and said to be the {\em link-cycle} or {\em link} of  $u $, denoted by $lk(u)$ or $lk_{\mathcal{K}}(u)$. 
A face in a map of the form  $u_1\mbox{-}u_2\mbox{-}\cdots\mbox{-}u_n\mbox{-}u_1$ is also denoted by  $[u_1,\cdots,u_n]$.

Let $ \mathcal{K}$ be a semi-equivelar map of type $X$.  We define a map $\mathcal{K}^{\ell}$ as follows. Let $ C $ be a non-trivial cycle (homotopic to a generator of the fundamental group and non-genus-separating) of length $ k $ in $ \mathcal{K} $ such that the cycle $C$ divides face-cycles of the vertices of $C$, that is, every sub-path $u\mbox{-}v\mbox{-}w \subset C$ of length two is a chord of the face-cycle($v$) at each vertex $v \in V(C)$. We cut $ \mathcal{K} $ along $ C $. We get a map $\mathcal{M}$ with two boundary cycles. Let $ A(p_1,\dots, p_k) $ and $ B(q_1,\dots, q_k) $ be the boundaries of $\mathcal{M}$ where $ p_i$ and $q_i$ ($ i=1,2,\dots, k $) are same in $\mathcal{K}$. Consider $ \mathcal{M}_1 $ and $ \mathcal{M}_2 $ where $ \mathcal{M}_1 \cong \mathcal{M}$, $ \mathcal{M}_2 \cong \mathcal{M}$. Then $ \mathcal{M}_i $ has two boundary cycles $ A_i(p_1^i,\dots, p_k^i) $ and $ B_i(q_1^i,\dots, q_k^i)$. We identify $ A_1(p_1^1,\dots, p_k^1) $ with $ B_2(q_1^2,\dots, q_k^2) $ by $ p_i^1 \rightarrow q_i^2 $, and $ B_1(q_1^1,\dots, q_k^1) $ with $ A_2(p_1^2,\dots, p_k^2) $ by $ q_i^1 \rightarrow p_i^2 $. Thus, we get a map $\mathcal{K}^2$ on the surface of \Echar{} $ 2\chi(\mathcal{K})$. Similarly, we consider $ \mathcal{M}_1, \dots, \mathcal{M}_{\ell}$ where $ \mathcal{M}_i \cong \mathcal{M}$ for each $i$. Let $ A_i, B_i$ denote the boundary cycles of $\mathcal{M}_i$. We identify $ B_i $ with $ A_{i+1} $ for $1 \le i \le {\ell}$ (addition in the suffix is modulo  ${\ell}$). Hence, we get a map $\mathcal{K}^{\ell}$ of type $X$, which is called a $ {\ell} $-covering map of $ \mathcal{K} $.  

Two maps of fixed type on a surface are {\em isomorphic} if there exists a {\em homeomorphism} of the surface, which maps vertices to vertices, edges to edges, faces to faces, and preserves incidents. More precisely, 
if we consider two polyhedral complexes $\mathcal{K}_{1}$ and $\mathcal{K}_{2}$ then an isomorphism to be a map $f ~:~ \mathcal{K}_{1}\rightarrow \mathcal{K}_{2}$ such that $f|_{V(\mathcal{K}_{1})} : V(\mathcal{K}_{1}) \rightarrow V(\mathcal{K}_{2})$ is a bijection and $f(\sigma)$ is a cell in $\mathcal{K}_{2}$ if and only if $\sigma$ is cell in $\mathcal{K}_{1}$. An isomorphism $f \colon \mathcal{K}\rightarrow \mathcal{K}$ is called an automorphism. The set $Aut(\mathcal{K})$ of all automorphism  of $\mathcal{K}$ forms a group {under the composition of map}, called the symmetric group of $\mathcal{K}$. A map $\mathcal{K}$ is said to be  {\em vertex-transitive} if {$Aut(\mathcal{K})$} acts transitively on the set  of vertices of $\mathcal{K}$. In \cite{lutz1999}, Lutz found all the (77 in numbers) vertex-transitive simplicial maps  with at most $15$ vertices. 

All vertex-transitive maps on the 2-sphere $\mathbb{S}^2$ are known. These are the boundaries of Platonic $\&$ Archimedean solids, the prism, and two infinite families (\cite{Ba1991}, \cite{GS1981}). Other than these, there exists a non-vertex-transitive semi-equivelar map on $\mathbb{S}^2$, namely the boundary of pseudorhombicuboctahedron (\cite{Gr2009}, \cite{wiki}).
There are eleven types of semi-equivelar maps on the torus, and all these are quotients of Archimedean tilings of the plane (\cite{DM2017}, \cite{datta2017}). 
Among these 11 types, four types ($[3^6]$, $[6^3]$, $[4^4]$, $[3^3, 4^2]$) of maps are always vertex-transitive, and there are infinitely many such examples in each type (\cite{Ba1991}, \cite{DM2017}). 
For each of the other seven types, there exists a semi-equivelar map on the torus, which is not vertex-transitive (\cite{DM2017}). Although, there are vertex-transitive maps of each of these seven types also (\cite{Ba1991}, \cite{Su2011t}, \cite{Th1991}). 
Similar results are known for Klein bottle (\cite{Ba1991}, \cite{DM2017}, \cite{Su2011kb}). 
If the Euler characteristic $\chi(M)$ of a surface $M$ is negative, then the number of semi-equivelar maps on $M$ is finite and at most $-84\chi(M)$ (\cite{Ba1991}).
Seventeen examples of non-vertex-transitive semi-equivelar maps 
on the surface with Euler characteristic $-1$ are known (\cite{DM2020}).  
There are precisely three non-vertex-transitive and seventeen vertex-transitive semi-equivelar maps on the orientable surface of genus 2 (\cite{DU2006, kn2012}), $103$ vertex-transitive semi-equivelar maps on the orientable surface of genus 3 (\cite{kn2012}), and 111 vertex-transitive semi-equivelar maps on the orientable surface of genus 4 (\cite{kn2012}). In this article, we prove the following.

\begin{theorem} \label{theo:sf2}
Let  $\mathcal{K}$ be a semi-equivelar map of type $X$ on the surface with Euler char.  $-2 $. Then, there exists a semi-equivelar map $\mathcal{K}^m$ of type $X $ on the surface with $\chi  = -2m$ for {$m = 2, 3, 4, \dots$}.
\end{theorem} 

\begin{theorem} \label{theo:sf3}
Let  $\mathcal{K}$ be a semi-equivelar map of type $X$ on the surface with Euler char.  $-2 $. Then, Aut$(\mathcal{K}^m) \cong \mathbb{D}_m$ (dihedral group of order $2m$) or $\mathbb{Z}_{m}$ (cyclic group of order $m$) for {$m = 2, 3, 4, \dots$}. 
\end{theorem}

\begin{remark}
\textnormal{Clearly, from Theorem \ref{theo:sf3}, $\mathcal{K}^m$ is non-vertex-transitive for {$m = 2, 3, 4, \dots$}.}
\end{remark}

It is a classical problem to ask: {what are} the types of semi-equivelar maps that exist on a surface. We know {the} complete list of types on a few surfaces, like, sphere, projective plane, torus, Klein bottle, and surface with $\chi = -1$. In this article, we discuss the same on the surface with Euler char. $-2$. From \cite{bu2019, DU2006, am2019, utm2014, kn2012, k2011}, we know the following. 

\begin{proposition}\label{prop1}
\justify
Let  $\mathcal{K}$ be a semi-equivelar map of type $X$ on the surface $S$ with Euler char.  $-2 $. If $X \in  \{[3^7]$, $[3^4, 4^2]$, $[4^1, 6^1, 16^1]$, $[4^1, 8^1, 12^1]$, $[6^2,  8^1]$, $[3^1,4^1,8^1,4^1]$, $[3^1, 6^1, 4^1, 6^1]$, $[4^3,6^1]$, $[3^4, 8^1]$, $[3^2, 4^1, 3^1, 6^1]$, $[3^1, 4^4]$, $[3^5, 4^1]$, $[3^1, 4^1, 3^1, 4^2]$, $[3^1, 5^3]$, $[3^3, 4^1, 3^1, 4^1]\}$
then there exists a map of type $X$ on $S$.
\end{proposition} 

Using Prop. \ref{prop1}, we prove the following.

\begin{theorem} \label{theo:sf0}
Let  $\mathcal{K}$ be a semi-equivelar map on the surface $S$ with Euler char.  $-2 $. Then, the type of $\mathcal{K}$ is 
\sloppypar{
 \hspace{-0.6cm}$[3^7]$, $[3^5,4^1]$, $[3^4,4^2]$, $[3^3,4^1,3^1, 4^1]$, $[3^5,5^1]$, $[3^4,7^1]$, $[3^4,8^1]$, $[3^4,9^1]$, $[3^4,10^1]$, $[3^2,5^1,3^1,5^1]$, 
 $[3^3,4^1,6^1]$, $[3^2,4^1,3^1,6^1]$,  $[3^1,4^1,3^1,4^2]$, $[3^1,4^4]$, $[3^1,7^1,3^1,7^1]$, $[3^1,4^1,7^1,4^1]$, $[3^1,4^1,8^1,4^1]$,
 $[4^3,5^1]$, $[3^1,4^1,9^1,4^1]$, $[3^1,4^1,10^1,4^1]$, $[3^1,5^3]$, $[3^1,6^1,4^1,6^1]$, $[4^1,6^1,8^1]$, $[4^3,6^1]$, $[5^1,4^1,5^1,4^1]$, 
 $[4^1,6^1,14^1]$, $[4^1,6^1,16^1]$, $[4^1,6^1,20^1]$, $[4^1,8^1,10^1]$, $[4^1,8^1,12^1]$, $[6^2,7^1]$, $[8^1, 6^2]$, $[9^1, 6^2]$, $[10^1, 6^2]$, $[5^1, 8^2]$, $[6^1, 8^2]$, $[4^1, 10^2]$, $[3^1, 14^2]$ or $[7^3]$.}
\end{theorem} 

\begin{cor}\label{cor1}
There exist a semi-equivelar map $\mathcal{K}^m$ of type
\sloppypar{
\hspace{-0.6cm}$[3^7]$, $[3^5,4^1]$, $[3^4,4^2]$, $[3^3,4^1,3^1, 4^1]$, $[3^5,5^1]$, $[3^4,7^1]$, $[3^4,8^1]$, $[3^4,9^1]$, $[3^4,10^1]$, $[3^2,5^1,3^1,5^1]$, 
$[3^3,4^1,6^1]$, $[3^2,4^1,3^1,6^1]$, $[3^1,4^1,3^1,4^2]$, $[3^1,4^4]$, $[3^1,7^1,3^1,7^1]$, $[3^1,4^1,7^1,4^1]$, $[3^1,4^1,8^1,4^1]$, $[4^3,5^1]$, $[3^1,4^1,9^1,4^1]$, $[3^1,4^1,10^1,4^1]$, $[3^1,5^3]$, $[3^1,6^1,4^1,6^1]$, $[4^1,6^1,8^1]$, $[4^3,6^1]$, $[5^1,4^1,5^1,4^1]$, $[4^1,6^1,14^1]$, $[4^1,6^1,16^1]$, $[4^1,6^1,20^1]$, $[4^1,8^1,10^1]$, $[4^1,8^1,12^1]$, $[6^2,7^1]$, $[8^1, 6^2]$, $[9^1, 6^2]$, $[10^1, 6^2]$, $[5^1, 8^2]$, $[6^1, 8^2]$, $[4^1, 10^2]$, $[3^1, 14^2]$  or $[7^3]$ on the surface of $\chi  = -2m$  and Aut$(\mathcal{K}^m) \cong \mathbb{D}_m$ or $\mathbb{Z}_{m}$ for {$m = 2, 3, 4, \dots$}.}
\end{cor}

The proof of Corollary \ref{cor1} follows from Theorems \ref{theo:sf2},  \ref{theo:sf3}, \ref{theo:sf0}.

\smallskip

In \cite{bu2019, bu2019_2, DU2006}, Bhowmik, Datta and Upadhyay used a well-known technique and classified (up to isomorphism) a few types of maps on $-1$ and $-2$ Euler characteristic's surfaces. The classification process is as follows:
First, it is considered a face-cycle of a vertex (for example, $v$) and then completed the face-cycles of all vertices appearing in the face-cycle of $v$ and continuing. By this process, one can classify every type of map listed in Theorem \ref{theo:sf0}. In this article, we have used slightly different arguments to classify the maps of type $[3^4,10^1]$ (see in the proof of Theorem \ref{theo:sf1}). This technique is efficient when $n_j=1$ for some $j$ in $[p_1^{n_1}, \dots, p_k^{n_k}] $, and the graph of the map can be drawn completely with all links of the vertices of all $p_j$-gons. Here, consider one $p_i$-gon ($i\neq j$) and check all the possibilities to identify $p_i$-gon with other {$p_i$-gons}.
Thus, we show the following. 

\begin{theorem} \label{theo:sf1}
Let {$\mathcal{K} $} be a semi-equivelar map of type  $[3^4, 10^1]$ on the surface $S$ with Euler char.  $-2 $. Then, $\mathcal{K} \cong \mathcal{K}_1[3^4,10^1], \mathcal{K}_2[3^4,10^1], \mathcal{K}_3[3^4,10^1]$ or $\mathcal{K}_4[3^4,10^1]$ (see Example \ref{eg:8maps-torus}).
\end{theorem}

\begin{remark}
\textnormal{One can classify other types of maps on the surface with Euler char.  $-2 $ by the {above} algorithms.
}
\end{remark}

\section{Examples}\label{eg:8maps-torus}
{\rm We present $25$ types of semi-equivelar maps on the surface {with} $\chi = -2 $.  $\mathcal{K}[4^1,5^1,4^1,5^1]$ {is represented by} the figure without red edges; $\mathcal{K}[3^2,5^1,3^1,5^1]$ {is represented by} the figure with red edges. $\mathcal{K}[3^1,4^1,9^1,4^1]$ {is represented by} the figure without red edges; $\mathcal{K}[3^4,9^1]$ {is represented by} the figure with red edges. }

\begin{figure}
    \centering
    \includegraphics[height= 7cm, width= 6cm]{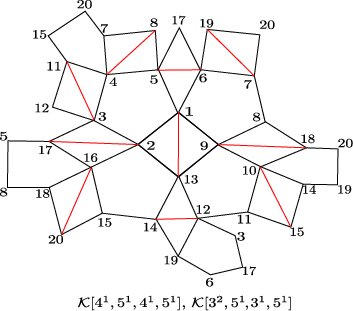}
    \includegraphics[height= 6cm, width= 8cm]{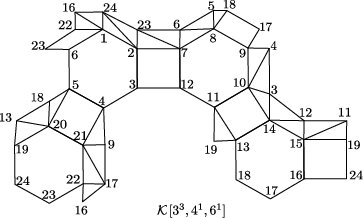}
    \label{4545_and_3(2)535}
\end{figure}
\hspace{-5cm}
\begin{figure}
    \centering
    \includegraphics[height= 6cm, width= 3cm]{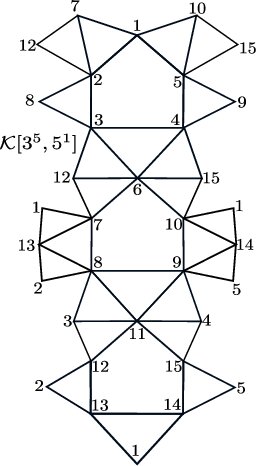}
    \includegraphics[height= 7cm, width= 11cm]{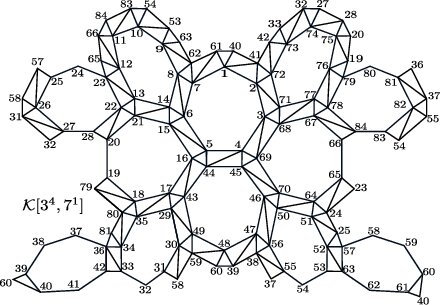}
\end{figure}
\hspace{-5cm}
\begin{figure}
    \centering
    \includegraphics[height= 5cm, width= 7cm]{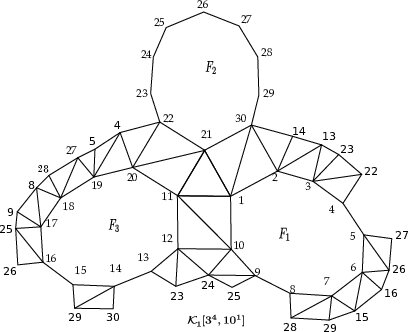}
    \includegraphics[height= 5cm, width= 7cm]{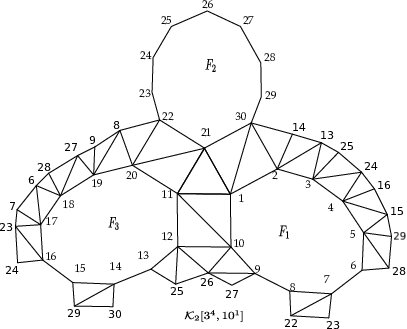}
\end{figure}

\begin{figure}
    \centering
    \includegraphics[height= 5cm, width= 7cm]{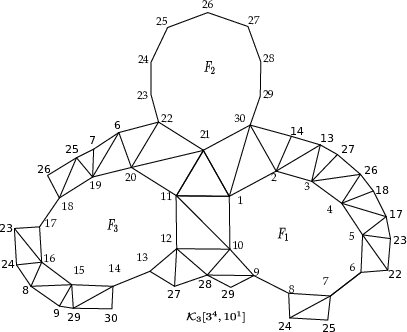}
    \includegraphics[height= 5cm, width= 7cm]{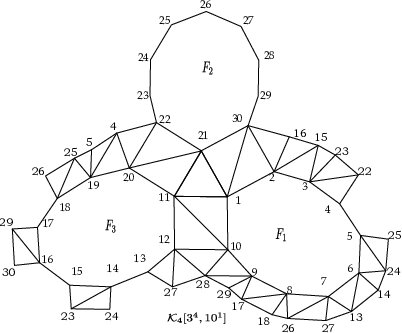}
\end{figure}

\begin{figure}
    \centering
    \includegraphics[height= 5cm, width= 6cm]{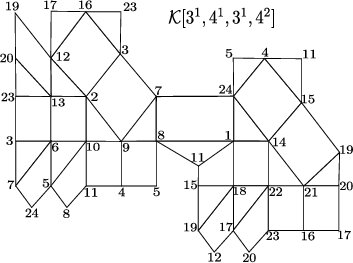}
    \includegraphics[height= 5cm, width= 9cm]{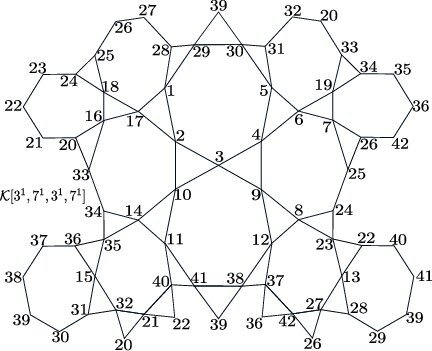}
\end{figure}

\begin{figure}
    \centering
    \includegraphics[height= 6cm, width= 11cm]{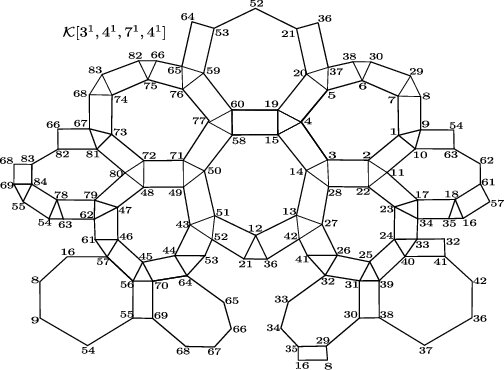}
\end{figure}

\begin{figure}
    \centering
    \includegraphics[height= 4cm, width= 6cm]{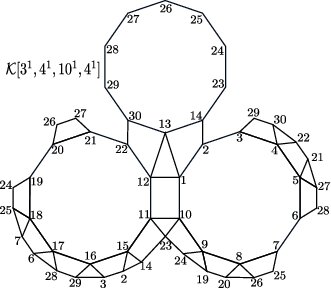}
    \includegraphics[height= 4cm, width= 8cm]{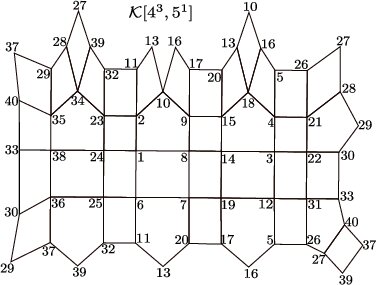}
\end{figure}

\begin{figure}
    \centering
    \includegraphics[height= 7cm, width= 12cm]{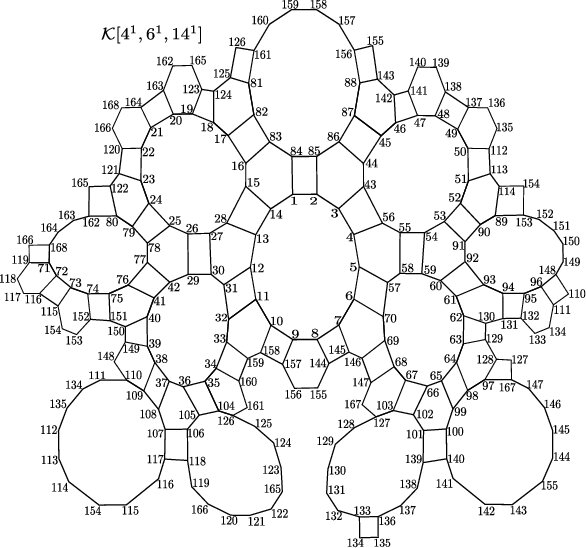}
\end{figure}

\begin{figure}
    \centering
    \includegraphics[height= 5cm, width= 7.5cm]{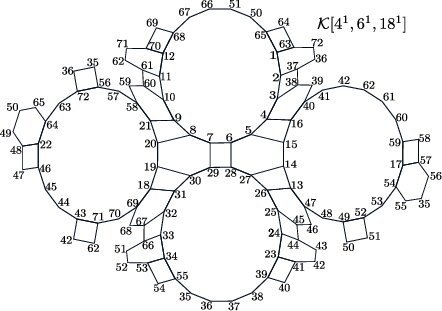}
    \includegraphics[height= 5cm, width= 7.5cm]{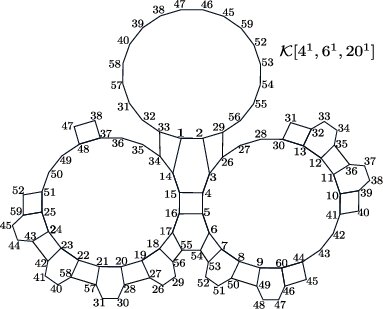}
\end{figure}

\begin{figure}
    \centering
    \includegraphics[height= 5cm, width= 8cm]{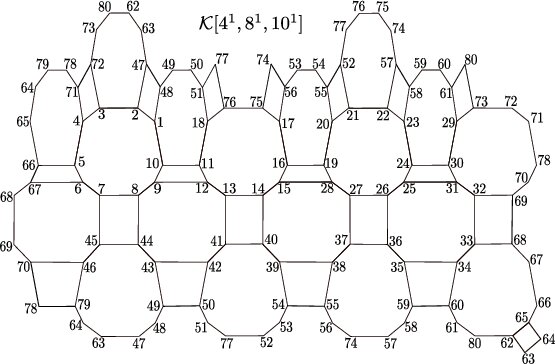}
    \includegraphics[height= 5cm, width= 7cm]{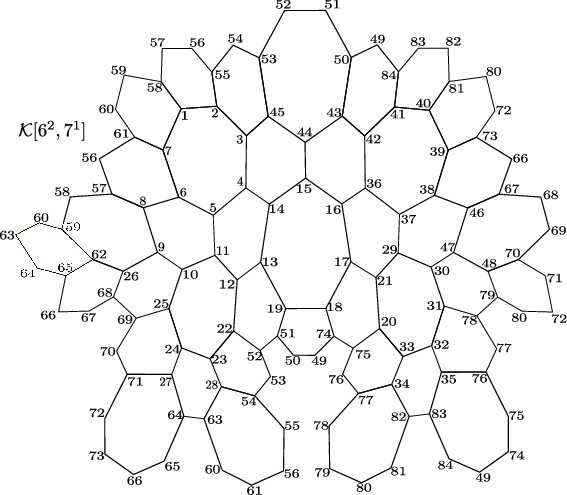}
\end{figure}

\begin{figure}
    \centering
    \includegraphics[height= 4cm, width= 6.5cm]{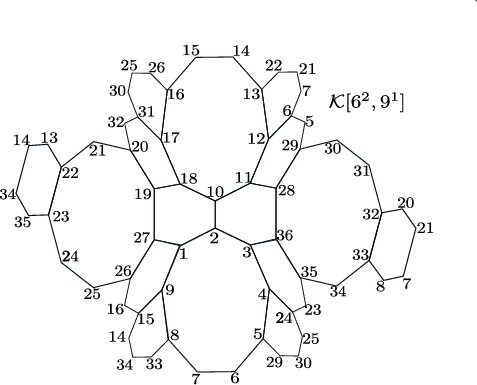}
    \includegraphics[height= 4cm, width= 6.5cm]{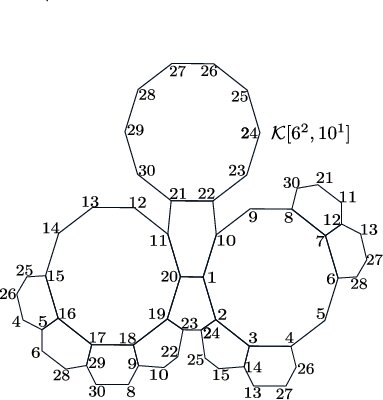}
\end{figure}

\begin{figure}
    \centering
    \includegraphics[height= 5cm, width= 8cm]{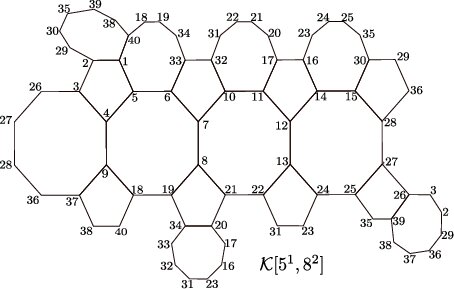}
    \includegraphics[height=6cm, width= 6cm]{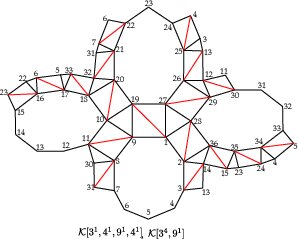}
\end{figure}
\begin{figure}
    \centering
    \includegraphics[height= 4.5cm, width= 7cm]{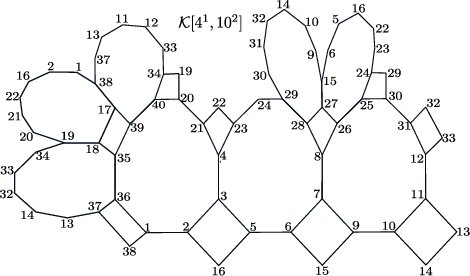}
    \includegraphics[height= 5cm, width= 6cm]{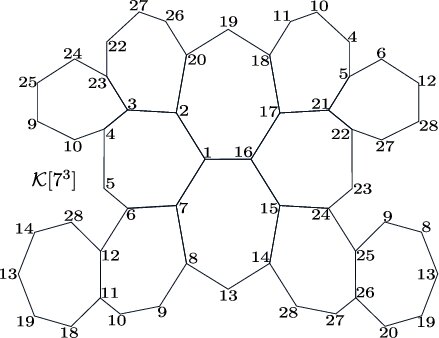}
\end{figure}
\begin{figure}
    \centering
    \includegraphics[height= 5cm, width= 5cm]{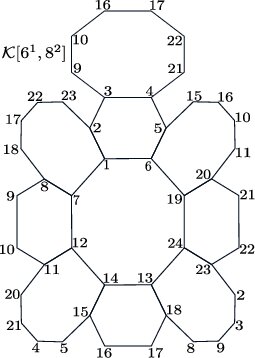}
    \includegraphics[height= 6cm, width= 10cm]{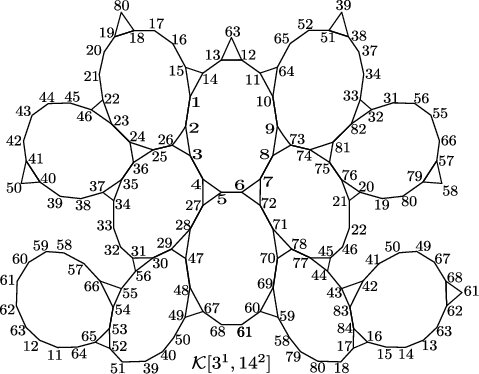}
\end{figure}

\begin{center}
\begin{longtable}{ |p{1cm}|p{2.2cm}|p{4.2cm}|p{2.4cm}|p{2.5cm}|}
\caption{Examples of semi-equivelar maps on the surface {with} $\chi = -2$}\\
\hline
S. No. & Type & Map &  Ref. & Orientation \\
\hline
\endfirsthead
\multicolumn{5}{c}%
{\tablename\ \thetable\ -- \textit{Continued from previous page}} \\
\hline
S. No. & Type & Map &  Ref. & Orientation \\
\hline
\endhead
\hline \multicolumn{5}{r}{\textit{Continued on next page}} \\
\endfoot
\hline
\endlastfoot
 \hline\hline
 1 & $[3^7] $ &  $ N_1, N_2, N_3, N_4, N_5, N_6 $ & \cite[Sec. 2]{DU2006} & Orientable\\
 \hline
 2 & $[3^5,4^1] $ &  $ A_{2.99} $  &  \cite[Appx. A]{kn2012} & Orientable \\
 \hline
 3 & $[3^4,4^2] $  &  $ A_{2.105}, A_{2.109}$  & \cite[Appx. A]{kn2012} & Orientable\\
 \hline
 4 & $[3^5,5^1] $  &  $\mathcal{K}[3^5,5^1]$  & & Non-orientable \\
 \hline
 5 & $[3^3,4^1,3^1,4^1]$ & $KO_{1x}, KO_{2x}, $  & \cite[Fig. 1]{bu2019} & Orientable,\\
  & ($=x$) & $ KNO_{1x}, KNO_{2x}, KNO_{3x}$  & & Non-orientable\\
\hline
6 & $[3^4,7^1]$ &  $\mathcal{K} [3^4,7^1]$  & & Orientable\\
\hline
7 & $[3^4,8^1]$ &  $ A_{2.89} $  & \cite[Appx. A]{kn2012} & Orientable\\
\hline
8 & $[3^4,9^1]$ &  $\mathcal{K}[3^4,9^1]$ &  & Orientable \\
\hline
9 & $[3^4,10^1] $ &  $\mathcal{K}_1[3^4,10^1], \mathcal{K}_2[3^4,10^1],\allowbreak \mathcal{K}_3[3^4,10^1], \mathcal{K}_4[3^4,10^1]$  & & Orientable \\
\hline
10 & $[3^2,5^1,3^1,5^1]$ &  $\mathcal{K}[3^2,5^1,3^1,5^1]$  & & Non-orientable \\
\hline
11 & $[3^3,4^1,6^1]$ &  $\mathcal{K}[3^3,4^1,6^1]$  & & Orientable \\
\hline
12 & $[3^2,4^1,3^1,6^1]$ &  $ A_{2.93} $  & \cite[Appx. A]{kn2012} & Orientable\\
\hline
13 & $[3^1,4^1,3^1,4^2]$ &  $\mathcal{K}[3^1,4^1,3^1,4^2]$  & & Non-orientable \\
\hline
14 & $[3^1,4^4]$ &  $ A_{2.97} $  & \cite[Appx. A]{kn2012} & Orientable\\
\hline
15 & $[3^1,7^1,3^1,7^1] $ & $\mathcal{K}[3^1,7^1,3^1,7^1]$   & & Orientable\\
\hline
16 & $[3^1,4^1,7^1,4^1]$ &  $\mathcal{K}[3^1,4^1,7^1,4^1]$  & & Non-orientable \\
\hline
17 & $[3^1,4^1,8^1,4^1]$ &  $ A_{2.33} $  & \cite[Appx. A]{kn2012} & Orientable\\
\hline
18 & $[3^1,4^1,9^1,4^1]$ &  $\mathcal{K}[3^1,4^1,9^1,4^1]$  & & Orientable \\
\hline
19 & $[3^1,4^1,10^1,4^1]$ &  $\mathcal{K}[3^1,4^1,10^1,4^1]$   & & Non-orientable\\
\hline
20 & $[3^1,5^3]$ &  (b), (c)  & \cite[Fig 4.4]{am2019} & Orientable\\
\hline
21 & $[3^1,6^1,4^1,6^1]$ &  $ A_{2.39} $  & \cite[Appx. A]{kn2012} & Orientable\\
\hline
22 & $[4^3,5^1]$ &  $\mathcal{K}[4^3,5^1]$  & & Non-orientable \\
\hline
23 & $[4^3,6^1]$ &  $ A_{2.43},  A_{2.46},  A_{2.51} $  & \cite[Appx. A]{kn2012} & Orientable\\
\hline
24 & $[5^1,4^1,5^1,4^1]$ &  $\mathcal{K}[4^1,5^1,4^1,5^1]$ & & Non-orientable \\
\hline
25 & $[4^1,6^1,14^1]$ &  $\mathcal{K}[4^1,6^1,14^1]$ & & Non-orientable \\
\hline
26 & $[4^1,6^1,16^1]$ &  $ A_{2.1} $ & \cite[Appx. A]{kn2012} & Orientable\\
\hline
27 & $[4^1,6^1,18^1]$ &  $ \mathcal{K}[4^1,6^1,18^1] $  & & Non-orientable \\
\hline
28 & $[4^1,6^1,20^1]$ &  $  \mathcal{K}[4^1,6^1,20^1] $ & & Non-orientable \\
\hline
29 & $[4^1,8^1,10^1]$ &  $\mathcal{K}[4^1,8^1,10^1]$ & & Non-orientable \\
\hline
30 & $[4^1,8^1,12^1]$ &  $ A_{2.7} $ & \cite[Appx. A]{kn2012} & Orientable\\
\hline
31 & $[7^1, 6^2]$ &  $\mathcal{K}[6^2,7^1]$ & & Non-orientable \\
\hline
32 & $[8^1, 6^2]$ &  $ A_{2.27} $ & \cite[Appx. A]{kn2012} & Orientable\\
\hline
33 & $[9^1, 6^2]$ & $ \mathcal{K}[6^2,9^1] $ & & Non-orientable \\
\hline
34 & $[10^1, 6^2]$ & $ \mathcal{K}[6^2,10^1] $ &  & Non-orientable \\
\hline
35 & $[5^1, 8^2]$ &  $\mathcal{K}[5^1, 8^2]$ & & Non-orientable \\
\hline
36 & $[6^1, 8^2]$ &  $\mathcal{K}[6^1,8^2]$ & & Non-orientable \\
\hline
37 & $[4^1, 10^2]$ &  $\mathcal{K}[4^1, 10^2]$ &  & Non-orientable \\
\hline
38 & $[3^1, 14^2]$ &  $\mathcal{K}[3^1,14^2]$  & & Orientable \\
\hline
39 & $[7^3] $ & $\mathcal{K}[7^3]$   & {\cite[Dual of $N_4$ (Sec. 2)]{DU2006}} & Orientable \\
\end{longtable}
\end{center}

\section{Semi-equivelar maps on the surface of $\chi = -2$} \label{sec:proofs-1}



We need the following technical proposition from \cite{DM2017} to prove Lemma  \ref{lem1}.

\begin{proposition} [Datta  $\& $ Maity] \label{prop}
If  $[p_1^{n_1}, \dots, p_k^{n_k}] $ satisfies any of the following three properties then  $[p_1^{n_1} $,  $\dots, p_k^{n_k}] $ can not be the type of any semi-equivelar map on a surface.
\begin{enumerate}[label=\roman*]
\item There exists  $i $ such that  $n_i=2 $,  $p_i $ is odd and  $p_j\neq p_i $ for all  $j\neq i $.

\item There exists  $i $ such that  $n_i=1 $,  $p_i $ is odd,  $p_j\neq p_i $ for all  $j\neq i $ and  $p_{i-1}\neq p_{i+1} $.

\item   $[p_1^{n_1} $,  $\dots, p_k^{n_k}] $ is of the form  $[p^1,q^m,p^1,r^n] $, where  $p $,  $q $,  $r $ are distinct and  $p $ is odd.
\end{enumerate}
(Here, addition in the subscripts are modulo  $k $.)
\end{proposition}

\begin{lemma} \label{lem1}
Let  $\mathcal{K} $ be an  $n $-vertex  map on the surface with Euler char.  $\chi = -2 $ of type   $[p_1^{n_1}, \dots,  $  $p_{\ell}^{n_{\ell}}] $.  Then,
 \begin{align*} 
 (n, [p_1^{n_1}, \dots, p_{\ell}^{n_{\ell}}]) \in & \{(12, [3^7]), (24,[3^5,4^1]),  (12,[3^4,4^2]), (12,[3^3,4^1,3^1,4^1]), (15,[3^5,5^1]), \\
 & (84,[3^4,7^1]), (48,[3^4,8^1]),(36,[3^4,9^1]), (30,[3^4,10^1]), (20,[3^2,5^1,3^1,5^1]),  \\
 & (40,[3^2,4^1,3^1,5^1]), (24,[3^3,4^1,6^1]), (24,[3^2,4^1,3^1,6^1]), (24,[3^1,4^1,3^1,4^2]),\\
 & (12,[3^1,4^4]),  (42,[3^1,7^1,3^1,7^1]), (24,[3^1,8^1,3^1,8^1]), (18,[3^1,9^1,3^1,9^1]), \\ 
 & (84,[3^1,4^1,7^1,4^1]), (48,[3^1,4^1,8^1,4^1]), (36,[3^1,4^1,9^1,4^1]), (30,[3^1,5^3]), \\
 &  (30,[3^1,4^1,10^1,4^1]),  (24,[3^1,6^1,4^1,6^1]), (15,[3^1,6^1,5^1,6^1]), (40,[4^3,5^1]), \\
 &  (24,[4^3,6^1]), (20,[5^1,4^1,5^1,4^1]),(168,[4^1,6^1,14^1]), (96,[4^1,6^1,16^1]), \\
 &  (72,[4^1,6^1,18^1]),(60,[4^1,6^1,20^1]),(80,[4^1,8^1,10^1]), (48,[4^1,8^1,12^1]),\\
 &   (84,[6^2, 7^1]),(48,[6^2, 8^1]), (36,[6^2, 9^1]), (30,[6^2,10^1]),(40,[5^1, 8^2]),\\
 &  (24,[6^1, 8^2]), (40,[4^1, 10^2]), (24,[4^1, 12^2]),(84,[3^1, 14^2]),\\
 &  (48,[3^1, 16^2]), (36,[3^1, 18^2]), (28,[7^3])\}.
 \end{align*}  
\end{lemma}

\begin{proof}
Let $f_0(=n)$, $f_1$, $f_2$ be the number of vertices, edges and faces of $\mathcal{K}$, respectively.
Let $d$ be the degree of each vertex. {Consider the $k$\mbox{-}tuple $(q_1^{m_1}, \dots, q_k^{m_k})$, where $3\le q_1 < \dots < q_k$, for each $i=1, \dots, k$, $q_i=p_j$ for some $j$, $m_i= \sum_{p_j=q_i} n_j$}. Let $x_i$ be the number of $i$-gons in $\mathcal{K}$. So, $\sum_i m_i = \sum_j n_j =d$ and $3f_2 \le 2f_1=f_0d$. Clearly, (the number of $i-$gons)  $\times  q_i = n \times m_i$, and hence $f_2=n\times(\frac{m_1}{q_1}+\dots +\frac{m_k}{q_k})$.
Since $\chi = -2$, we get, $-2-f_0 = f_2-f_1\le (f_0 \times d)/3-f_1=-(f_0\times d)/6$. Thus, $(d-6)f_0\le -6 \chi$ and hence $d\le 7$ since $f_0 \ge 12$. So, $d= 3,4,5,6$ or $7$.

\smallskip

\noindent \textbf{Case 1.} Assume $d=7$. Then, $ \sum n_i=7$. Therefore, $7f_0=2f_1=\sum (i\times x_i)$, and $f_2=\sum x_i$.
Now, $-2=f_0-f_1+f_2=(2/7)f_1-f_1+f_2 = (-5/7)f_1+f_2$, therefore, $14=5f_1-7f_2$, i.e., $28=10f_1-14f_2=\sum (5i-14)x_i=x_3+6x_4$ since $x_i\ge3$ for all $i$. Therefore, $(x_3,x_4)=(28,0), (10,3), (3,4).$\\
Let $(x_3,x_4)=(28,0)$, {then} $(q_1^{m_1}, \dots, q_k^{m_k})=(3^{l_1})$. {Therefore, $l_1=7$ and hence $[p_1^{n_1}, \dots, p_{\ell}^{n_{\ell}}]\allowbreak=[3^7]$ with $f_0=12$.}\\
Let $(x_3,x_4)=(10,3)${, then} $(q_1^{m_1}, \dots, q_k^{m_k})=(3^{l_1},4^{l_2})$, where $l_1+l_2=7, l_1 >0, l_2 >0$.
{Therefore,} {$(l_1f_0)/3=10$} and {$(l_2f_0)/4=3$}.
Hence, $-2=f_0-f_1+f_2=f_0-(7f_0)/2+(l_1/3+l_2/4)f_0$, i.e., {$f_0=6$, which is not possible}. 
Hence, $(x_3,x_4) \neq (10,3)$. \\
Similarly, $(x_3,x_4) \neq (3,4)$.


\smallskip

\noindent \textbf{Case 2.}
Assume $d=6$. Then, $ \sum n_i=6$. Therefore, $6f_0=2f_1=\sum (i\times x_i)$ and $-2=f_0-f_1+f_2=f_1/3-f_1+f_2=(-2/3)f_1+f_2$, therefore, $6=2f_1-3f_2=\sum (i\times x_i)-3\sum x_i=\sum (i-3)x_i$, i.e., $6=x_4+2x_5+3x_6+4x_7+5x_8+6x_9+\dots$. Hence, $x_i=0$ for all $i\ge 7$. 
So, $(x_4,x_5,x_6,x_7,x_8,x_9)=(6,0,0,0,0,0),(0,3,0,0,0,0)$, since {$x_i\ge 3$} for all $i$.\\
If $(x_4,x_5,x_6,x_7,x_8,x_9)=(6,0,0,0,0,0)$, then $(q_1^{m_1}, \dots, q_k^{m_k})=(3^{l_1},4^{l_2})$ and $l_1+l_2=6$. So, $-2=n-6n/2+nl_1/3+nl_2/4=-2n+n(l_1/3+l_2/4)$, i.e., $24=n(6-l_1)$. Thus $(l_1,l_2,n)=(5,1,24),(4,2,12)$ since $n\ge 12$. Therefore,  $(q_1^{m_1}, \dots, q_k^{m_k})=(3^5,4^1),(3^4,4^2)$. So, $[p_1^{n_1}, \dots, p_{\ell}^{n_{\ell}}]=[3^5,4^1],[3^4,4^2],[3^3,4^1,3^1,4^1],[3^2,4^1,3^2,4^1]$. However, $[3^2,4^1,3^2,4^1]$ is not possible by Prop 3.1. Hence, $[p_1^{n_1}, \dots, p_{\ell}^{n_{\ell}}]=[3^5,4^1],[3^4,4^2],[3^3,4^1,3^1,4^1]$ \\
If $(x_4,x_5,x_6,x_7,x_8,x_9)=(0,3,0,0,0,0)$, then $(q_1^{m_1}, \dots, q_k^{m_k})=(3^{l_1},5^{l_2})$ and $l_1+l_2=6$. So, $-2=n-6n/2+nl_1/3+nl_2/5=-2n+n(l_1/3+l_2/5)$, i.e., $15=n(6-l_1)$. Thus $(l_1,l_2,n)=(5,1,15)$ since $n\ge 12$. Therefore, $(q_1^{m_1}, \dots, q_k^{m_k})=(3^5,5^1)$. Hence, $[p_1^{n_1}, \dots, p_{\ell}^{n_{\ell}}]=[3^5,5^1]$

\smallskip
\noindent \textbf{Case 3.}
Assume $d=5$. Then, $ \sum n_i=5$. Therefore, $5f_0=2f_1=\sum (i\times x_i)$ and $-2=f_0-f_1+f_2=(2/5)f_1-f_1+f_2=(-3/5)f_1+f_2$, therefore, $10=3f_1-5f_2$, i.e., $20=6f_1-10f_2=\sum (3i-10)x_i=\sum_{i\ge 4} (3i-10)x_i-x_3$\\
Let $(q_1^{m_1}, \dots, q_k^{m_k})=(3^4,i^1)$. Then, $x_3=(4/3)f_0$, $x_i=f_0/i$. Hence, $20=(3i-10)f_0/i-(4/3)f_0$, i.e., $f_0=12i/(i-6)$.
Therefore, $(i,f_0,x_3,x_i)=(7,84,112,12), (8,48,64,6), $ $(9,36,48,4), (10,30,40,3)$ as $x_i\ge 3$ for all $i$, and $f_0 \ge 12$. So,  $(q_1^{m_1}, \dots, q_k^{m_k})= (3^4,7^1), \allowbreak (3^4,8^1), \allowbreak (3^4,9^1), \allowbreak (3^4,10^1)$. Hence, $[p_1^{n_1}, \dots, p_{\ell}^{n_{\ell}}]=[3^4,7^1], [3^4,8^1], [3^4,9^1], [3^4,10^1]$.
\\
Let $(q_1^{m_1}, \dots, q_k^{m_k})=(3^3,i^2)$. Then, $x_3=f_0$, $x_i=2f_0/i$, hence, $20=(3i-10)\times 2f_0/i-f_0$, i.e., $f_0=4i/(i-4)$.
Therefore, $(i,f_0,x_3,x_i)=(5,20,20,8), (6,12,12,4)$ as $x_i\ge 3$ for all $i$, and $f_0 \ge 12$.  So,  $(q_1^{m_1}, \dots, q_k^{m_k})= (3^3,5^2), (3^3,6^2)$. Hence, $[p_1^{n_1}, \dots, p_{\ell}^{n_{\ell}}]= [3^3,5^2],[3^2,5^1,3^1,5^1],[3^3,6^2], \allowbreak [3^2,6^1,3^1,6^1]$. However, $[3^3,5^2]$ is not possible by Prop 3.1. Also, $[3^3,6^2]$ and $[3^2,6^1,3^1,6^1]$ are not possible as the edge graphs of are not complete. Therefore, we have, $[p_1^{n_1}, \dots, p_{\ell}^{n_{\ell}}]= [3^2,5^1,3^1,5^1]$.
\\
Let $(q_1^{m_1}, \dots, q_k^{m_k})=(3^3,i^1,j^1)$, $i<j$. Then, $x_3=f_0$, $x_i=f_0/i$, $x_j=f_0/j$, hence, $20=(3i-10)f_0/i+(3j-10)f_0/j-f_0$, i.e., $f_0=4ij/(ij-2i-2j)$.
Therefore, $(i,j,f_0,x_3,x_i,x_j)=(4,5,40,40,10,8), (4,6,24,24,6,4)$ as $x_i\ge 3$ for all $i$, and $f_0 \ge 12$.  So,  $(q_1^{m_1}, \dots, q_k^{m_k})= (3^3,4^1,5^1), \allowbreak (3^3,4^1,6^1)$. Hence, $[p_1^{n_1}, \dots, p_{\ell}^{n_{\ell}}]=[3^3,4^1,5^1]$, $  [3^2,4^1,3^1,5^1], [3^3,4^1,6^1], \allowbreak [3^2,4^1,3^1,\allowbreak 6^1]$. However, $[3^3,4^1,5^1]$ is not possible by Prop 3.1. Therefore, we have, $[p_1^{n_1}, \dots, p_{\ell}^{n_{\ell}}]= [3^2,4^1,3^1,5^1], \allowbreak [3^3,4^1,6^1], [3^2,4^1,3^1,6^1]$.
\\
Let $(q_1^{m_1}, \dots, q_k^{m_k})=(3^2,i^1,j^1,a^1)$, $i<j<a$. By Prop. 3, $[p_1^{n_1}, \dots, p_{\ell}^{n_{\ell}}]= [3^1,i^1,3^1,j^1,i^1],$ $[3^1,i^1,3^1,i^2]$, where $i>3,$ $j>3$.\\
Let $[p_1^{n_1}, \dots, p_{\ell}^{n_{\ell}}]= [3^1,i^1,3^1,j^1,i^1]$. Then, $x_3=(2/3)f_0$, $x_i=2f_0/i$, $x_j=f_0/j$, hence, $20=(3i-10)\times 2f_0/i+(3j-10)f_0/j-(2/3)f_0$, i.e., $f_0=12ij/(5ij-6i-12j)$.
Therefore, $(i,j,f_0,x_3,x_i,x_j)=(4,4,24,16,12,6), (6,3,12,8,4,4)$ as $x_i\ge 3$ for all $i$, and $f_0 \ge 12$. So, $[p_1^{n_1}, \dots, p_{\ell}^{n_{\ell}}]=[3^1,4^1,3^1,4^2], [3^1,6^1,3^2,6^1]$. However, $[3^1,6^1,3^2,6^1]$ is not possible as mentioned above. Therefore, we have, $[p_1^{n_1}, \dots, p_{\ell}^{n_{\ell}}]= [3^1,4^1,3^1,4^2]$.\\
Let $[p_1^{n_1}, \dots, p_{\ell}^{n_{\ell}}]= [3^1,i^1,3^1,i^2]$. Then, $x_3=(2/3)f_0$, $x_i=3f_0/i$, hence, $20=(3i-10)\times 3f_0/i-(2/3)f_0$, i.e., $f_0=12i/(5i-18)$.
Therefore, $(i,f_0,x_3,x_i)=(4,24,16,18)$ as $x_i\ge 3$ for all $i$, and $f_0 \ge 12$. So, $[p_1^{n_1}, \dots, p_{\ell}^{n_{\ell}}]= [3^1,4^1,3^1,4^2]$.
\\
Let $(q_1^{m_1}, \dots, q_k^{m_k})=(3^1,i^1,j^1,a^1,b^1)$, $i<j<a<b$. By Prop. 3, $[p_1^{n_1}, \dots, p_{\ell}^{n_{\ell}}]= [3^1,i^1,j^1,a^1,i^1], [3^1,i^1,j^1,i^2]$, where $i>3$, $j>3$, $a>3$. \\
Let $[p_1^{n_1}, \dots, p_{\ell}^{n_{\ell}}]= [3^1,i^1,j^1,a^1,i^1]$. Then, $x_3=f_0/3$, $x_i=2f_0/i$, $x_j=f_0/j$, $x_a=f_0/a$, hence, $20=(3i-10)\times 2f_0/i+(3j-10)f_0/j+(3a-10)f_0/a-f_0/3$, i.e., $f_0=12ija/(7ija-6ia-6ij-12ja)$.
Therefore, $(i,j,a,f_0,x_3,x_i,x_j,x_a)=(4,4,3,24,8,12,6,8),$ $(4,4,4,12,4,6,3,3)$ as $x_i\ge 3$ for all $i$, and $f_0 \ge 12$. So, $[p_1^{n_1}, \dots, p_{\ell}^{n_{\ell}}]= [3^1,4^2,3^1,4^1],[3^1,4^4]$.

\smallskip

\noindent \textbf{Case 4.}
Assume $d=4$. Then, $ \sum n_i=4$. Therefore, $4f_0=2f_1=\sum (i\times x_i)$ and $-2=f_0-f_1+f_2=f_1/2-f_1+f_2 = (-1/2)f_1 + f_2$, therefore, $4=f_1-2f_2$, i.e., $8=2f_1-4f_2=\sum (i-4)x_i=\sum_{i\ge 4} (i-4)x_i-x_3$\\
Let $(q_1^{m_1}, \dots, q_k^{m_k})=(3^3,i^1)$. Then, $x_3=f_0$, $x_i=f_0/i$. hence, $8=(i-4)f_0/i-f_0$, i.e., $f_0=-2i$, which is not possible. 
\\
Let $(q_1^{m_1}, \dots, q_k^{m_k})=(3^2,i^1,j^1)$. By Prop. 3, $[p_1^{n_1}, \dots, p_{\ell}^{n_{\ell}}]= [3^1,i^1,3^1,i^1],$, where $i>3$. 
Then, $x_3=(2/3)f_0$, $x_i=2f_0/i$, hence, $8=(i-4)\times 2f_0/i-(2/3)\times f_0$, i.e., $f_0=6i/(i-6)$.
Therefore, $(i,f_0,x_3,x_i)=(7,42,28,12),(8,24,16,6),(9,18,12,4),(10,15,10,13)$ as $x_i\ge 3$ for all $i$, and $f_0 \ge 12$. So, $[p_1^{n_1}, \dots, p_{\ell}^{n_{\ell}}]= [3^1,7^1,3^1,7^1],[3^1,8^1,3^1,8^1],[3^1,9^1,3^1,9^1],[3^1,10^1,\allowbreak 3^1,10^1]$. However, $[3^1,10^1,3^1,10^1]$ is not possible as the edge graph is not complete. Therefore, we have, $[p_1^{n_1}, \dots, p_{\ell}^{n_{\ell}}]= [3^1,7^1,3^1,7^1],[3^1,8^1,3^1,8^1],[3^1,9^1,3^1,9^1]$.
\\
Let $(q_1^{m_1}, \dots, q_k^{m_k})=(3^1,i^1,j^1,a^1)$. By Prop. 3, $[p_1^{n_1}, \dots, p_{\ell}^{n_{\ell}}]= [3^1,i^1,j^1,i^1],$, where $i>3$ and $j>3$.
Then, $x_3=f_0/3$, $x_i=2f_0/i$, $x_j=f_0/j$, hence, $8=(i-4)\times 2f_0/i + (j-4)\times f_0/j-f_0/3$, i.e., $f_0=6ij/(2ij-3i-6j)$. 
Therefore, $(i,j,f_0,x_3,x_i,x_j)= $ $(4,7,84,28,42,12), $ $(4,8,48,16,24,6), (4,9,36,12,18,4), (4,10,30,10,15,3), (5,4,120,40,48,30),$ $(5,5,30,10,12,\allowbreak 6),$ $(6,4,24,8,8,6), (6,5,15,5,5,3), (8,4,12,4,3,3)$ as $x_i\ge 3$ for all $i$, and $f_0 \ge 12$. 
So, $[p_1^{n_1}, \dots, p_{\ell}^{n_{\ell}}]= [3^1,4^1,7^1,4^1],$ $[3^1,4^1,8^1,4^1], [3^1,4^1,9^1,4^1],$ $[3^1,4^1,10^1,4^1], [3^1,5^1,4^1,5^1],$ $[3^1,\allowbreak 5^3],$ $[3^1,6^1,4^1,\allowbreak 6^1],$ $[3^1,6^1,5^1,6^1],$ $[3^1,8^1,4^1,8^1]$.
However, $[3^1,5^1,4^1,5^1]$ is not possible by Prop 3.1., and  $[3^1,8^1,4^1,8^1]$ is not possible since face-cycle of a vertex contains 16 vertices and $16>12 (=f_0)$.
Therefore, we have, $[p_1^{n_1}, \dots, p_{\ell}^{n_{\ell}}]= [3^1,4^1,7^1,4^1],$ $[3^1,4^1,9^1,4^1],$ $[3^1,4^1,10^1,4^1], [3^1,4^1,\allowbreak 8^1,4^1],$ $[3^1,5^3],$ $[3^1,6^1,4^1,6^1],$ $[3^1,6^1,5^1,6^1]$.
\\
Let $(q_1^{m_1}, \dots, q_k^{m_k})=(i^1,j^1,a^1,b^1)$, where $3<i<j<a<b$. Then, $x_i=f_0/i$, $x_j=f_0/j$, $x_a=f_0/a$, $x_b=f_0/b$, hence, $8=(i-4)f_0/i + (j-4)f_0/j + (a-4)f_0/a + (b-4)f_0/b$, i.e., $f_0=2ijab/(ijab-jab-iab-ijb-ija)$.
However, $(q_1^{m_1}, \dots, q_k^{m_k})=(i^1,j^1,a^1,b^1)$ is not possible for $3<i<j<a<b$, since $f_0\ge 12$.
\\
Let $(q_1^{m_1}, \dots, q_k^{m_k})=(i^2,j^1,a^1)$, where $i>3$, $j>3$ and $a>3$.
If $i$ is even, by Prop. 3, $[p_1^{n_1}, \dots, p_{\ell}^{n_{\ell}}]= [i^2,j^1,a^1], [i^1,j^1,i^1,a^1]$, where $i>3$, $j>3$ and $a>3$.
If $i$ is odd, by Prop. 3, $[p_1^{n_1}, \dots, p_{\ell}^{n_{\ell}}]= [i^1,j^1,i^1,j^1]$, where $i>3$ and $j>3$.
\\
Let $[p_1^{n_1}, \dots, p_{\ell}^{n_{\ell}}]= [i^2,j^1,a^1]$, where $i$ is even. Then, $x_i=2f_0/i$, $x_j=f_0/j$, $x_a=f_0/a$, hence, $8=(i-4)\times 2f_0/i + (j-4)\times f_0/j + (a-4)\times f_0/a$, i.e., $f_0=2ija/(ija-ij-ia-2ja)$.
Therefore, $(i,j,a,f_0,x_i,x_j,x_a)= $ $(4,4,5,40,20,10,8), $ $(4,4,6,24,12,6,4),$ $(4,5,5,20,10,4,4), (5,4,4,20,8,5,4), $ $(6,4,4,12,4,3,3)$, as $x_i\ge 3$ for all $i$, and $f_0 \ge 12$..
So, $[p_1^{n_1}, \dots, p_{\ell}^{n_{\ell}}]= [4^3,5^1], [4^3,6^1], \allowbreak [4^2,5^2], [6^2,4^2]$.
However, $[4^2,5^2]$ is not possible by Prop 3.1., and  $[6^2,4^2]$ is not possible since face-cycle of a vertex contains 13 vertices and $13>12 (=f_0)$.
Therefore, we have, $[p_1^{n_1}, \dots, p_{\ell}^{n_{\ell}}]= [4^3,5^1], [4^3,6^1].$
\\
Let $[p_1^{n_1}, \dots, p_{\ell}^{n_{\ell}}]= [i^1,j^1,i^1,a^1]$, where $i$ is even. Then, $x_i=2f_0/i$, $x_j=f_0/j$, $x_a=f_0/a$, hence, $8=(i-4)\times 2f_0/i + (j-4)\times f_0/j + (a-4)\times f_0/a$, i.e., $f_0=2ija/(ija-ij-ia-2ja)$.
Therefore, $(i,j,a,f_0,x_i,x_j,x_a) = (4,4,5,40,10,10,8), (4,4,6,24,6,6,4)$, as $x_i\ge 3$ for all $i$, and $f_0 \ge 12$.. 
So, $[p_1^{n_1}, \dots, p_{\ell}^{n_{\ell}}]= [4^3,5^1], [4^3,6^1]$.
\\
Let $[p_1^{n_1}, \dots, p_{\ell}^{n_{\ell}}]= [i^1,j^1,i^1,j^1]$, where $i$ is odd.
Then, $x_i=2f_0/i$, $x_j=2f_0/j$, hence, $8=(i-4)\times 2f_0/i + (j-4)\times 2f_0/j$, i.e., $f_0=2ij/(ij-2i-2j)$.
Therefore, $(i,j,f_0,x_i,x_j)= (5,4,20,8,10)$, as $x_i\ge 3$ for all $i$, and $f_0 \ge 12$.. 
So, $[p_1^{n_1}, \dots, p_{\ell}^{n_{\ell}}]= [5^1,4^1,5^1,4^1]$.
\\
Let $(q_1^{m_1}, \dots, q_k^{m_k})=(i^3,j^1)$, where $i>3$ and  $j>i$.
Then, $x_i=3f_0/i$, $x_j=f_0/j$, hence, $8=(i-4)\times 3f_0/i + (j-4)\times f_0/j$, i.e., $f_0=2ij/(ij-i-3j)$.
Therefore, $(i,j,f_0,x_i,x_j)= (4,5,40,10,8), (4,6,24,6,4)$, as $x_i\ge 3$ for all $i$, and $f_0 \ge 12$.
So,  $(q_1^{m_1}, \dots, q_k^{m_k})=(4^3,5^1), (4^3,6^1)$
Hence, $[p_1^{n_1}, \dots, p_{\ell}^{n_{\ell}}]= [4^3,5^1], [4^3,6^1]$.

\smallskip

\noindent \textbf{Case 5.}
Assume $d=3$. Then, $ \sum n_i=3$. Therefore, $3f_0=2f_1=\sum (i\times x_i)$ and $-2=f_0-f_1+f_2=(2/3)f_1-f_1+f_2 =(-1/3)f_1+f_2$, therefore, $6=f_1-3f_2$, i.e., $12=2f_1-6f_2=\sum (i-6)x_i$\\
Let $(q_1^{m_1}, \dots, q_k^{m_k})=(i^1,j^1,a^1)$, where $i<j<a$. By Prop. 3, $[p_1^{n_1}, \dots, p_{\ell}^{n_{\ell}}] = [i^1,j^1,a^1]$, where $i<j<a$, and $i,j$ and $a$ are even. Then, $x_i=f_0/i$, $x_j=f_0/j$ and $x_a=f_0/a$. hence, $12=(i-6)f_0/i+(j-6)f_0/j+(a-6)f_0/a$, i.e., $f_0=4ija/(ija-2ij-2ia-2ja)$.
Therefore, $(i,j,a,f_0,x_i,x_j,x_a) =$ $ (4,6,14,168,42,28,12),$ $(4,6,16,96,24,16,6),$ $ (4,6,18,72,18,12,4), $ $(4,6,20,60,15,10,3), $ $(4,8,10,80,20,10,8),$ $ (4,8,12,48,12,6,4).$
So, $[p_1^{n_1}, \dots, p_{\ell}^{n_{\ell}}]= [4^1,6^1,14^1], [4^1,6^1,16^1], [4^1,6^1,18^1], [4^1,6^1,20^1], [4^1,8^1,10^1], [4^1,8^1,12^1].$
\\
Let $(q_1^{m_1}, \dots, q_k^{m_k})=(i^2,j^1)$.
If $i$ is even, by Prop. 3, $[p_1^{n_1}, \dots, p_{\ell}^{n_{\ell}}]= [i^2,j^1]$.
If $i$ is odd, by Prop. 3, $[p_1^{n_1}, \dots, p_{\ell}^{n_{\ell}}]= [i^3]$.\\
Let $[p_1^{n_1}, \dots, p_{\ell}^{n_{\ell}}]= [i^2,j^1]$, where $i$ is even. 
Then, $x_i=2f_0/i$ and $x_j=f_0/j$, hence, $12=(i-6)\times 2f_0/i+(j-6)\times f_0/j$, i.e., $f_0=4ij/(ij-2i-4j)$.
Therefore, $(i,j,f_0,x_i,x_j)= (6,7,84,28,12), (6,8,48,6,16), (6,9,36,12,4), (6,10,30,10,3),$ $(8,5,40,10,8),$ $ (8,6,24,6,4),$ $(10,4,40,8,10),$ $(10,5,20,4,4),$ $ (12,4,24,4,6), $ $(14,3,84,12,28),$ $ (16,3,48,6,18)$  $(18,3,36,\allowbreak 4,12),$ $(20,3,30,3,10).$
So, $[p_1^{n_1}, \dots, p_{\ell}^{n_{\ell}}]= [6^2,7^1], [6^2,8^1], [6^2,9^1], [6^2,10^1], [8^2,5^1], [8^2,6^1],$ $ [10^2,4^1],$ $ [10^2,5^1], [12^2,4^1], [14^2,3^1], [16^2,3^1], [18^2,3^1], [20^2,3^1].$
However, $ [10^2,5^1]$ is not possible as the edge graph of $\mathcal{K}$ is not complete. Also,  $[20^2,3^1]$ is not possible as the face-cycle of a vertex contains $38$ vertices and $38>30 (=f_0)$.
So, $[p_1^{n_1}, \dots, p_{\ell}^{n_{\ell}}]= [6^2,7^1], [6^2,8^1], [6^2,9^1],$  $[6^2, 10^1], $ $[8^2,5^1], [8^2,6^1], [10^2,4^1],$ $ [12^2,4^1], [14^2,3^1], [16^2,3^1], [18^2,3^1].$
\\
Let $[p_1^{n_1}, \dots, p_{\ell}^{n_{\ell}}]= [i^3]$, where $i$ is odd.
Then, $x_i=3f_0/i$, hence, $12=(i-6)\times 3f_0/i$, i.e., $f_0=4i/(i-6)$.
Therefore, $(i,f_0,x_i)= (7,28,12), (9,12,4)$. So, $[p_1^{n_1}, \dots, p_{\ell}^{n_{\ell}}]= [7^3],[9^3]$.
However, $[9^3]$ is not possible since face-cycle of a vertex contains 21 vertices and $21>12 (=f_0)$. Therefore, we have, $[p_1^{n_1}, \dots, p_{\ell}^{n_{\ell}}]= [7^3]$. 

It completes the proof.
\end{proof}

\begin{lemma}\label{lem2}
	{Let  $\mathcal{K}$ be a SEM of type  $X$ on the surface of Euler char. $-2$. Then,  $X \neq [3^1, 8^1, 3^1, 8^1], [3^1, 9^1, 3^1, 9^1], [3^1, 6^1, 5^1, 6^1],  [4^1, 12^2], [3^1, 16^2], [3^1, 18^2]$.}
\end{lemma}

\begin{proof}
Let  $\mathcal{K}$ be a SEM of type  $[3^1, 8^1, 3^1, 8^1]$ on the surface of Euler char. $-2$. Let $u$ be a vertex in $\mathcal{K}$. Let $F_{1, 3} = uvw, F_{2, 3}=uxy$ be two  $3$-gons at $u$. Clearly, there are exactly four $8$-gons adjacent to $F_{1, 3} \cup F_{2,3}$, namely, $F_{1, 8}, F_{2, 8}, F_{3, 8}, F_{4, 8}$. Let $uv \in F_{1, 8}, vw \in F_{2, 8}, wu \in F_{3, 8}, uy \in F_{3, 8}, yx \in F_{4, 8}, xu \in F_{1, 8}$. Clearly, $F_{1, 8} \cap F_{2, 8} = v, F_{2, 8} \cap F_{3, 8} = w, F_{3, 8} \cap F_{4, 8} = y, F_{4, 8} \cap F_{1, 8} = x, F_{1, 8} \cap F_{3, 8} = u, \# V(F_{2, 8}) \cap V(F_{4, 8}) \le 1$, where $\# V(F_{2, 8}) \cap V(F_{4, 8})$ denotes the cardinality of the set $V(F_{2, 8}) \cap V(F_{4, 8})$. Hence, $\#V(\mathcal{K}) \ge 8+7+6+5 = 26$. This is a contradiction as   $\#V(\mathcal{K}) =24$. Hence,  $X \neq [3^1, 8^1, 3^1, 8^1]$.

Similarly, $X \neq [3^1, 9^1, 3^1, 9^1], [3^1, 6^1, 5^1, 6^1],  [4^1, 12^2], [3^1, 16^2], [3^1, 18^2]$.
\end{proof}

\begin{proof}[Proof of Theorem \ref{theo:sf0}]
Let $\mathcal{K}$ be a SEM on the surface {with} Euler char. $-2$ of type  $X$.  Then,  by Lemma \ref{lem1}, 
\begin{align*} 
(n, X) \in & \{(12, [3^7]), (24,[3^5,4^1]),  (12,[3^4,4^2]), (12,[3^3,4^1,3^1,4^1]), (15,[3^5,5^1]), \\
& (84,[3^4,7^1]), (48,[3^4,8^1]),(36,[3^4,9^1]), (30,[3^4,10^1]), (20,[3^2,5^1,3^1,5^1]),  \\
& (40,[3^2,4^1,3^1,5^1]), (24,[3^3,4^1,6^1]), (24,[3^2,4^1,3^1,6^1]), (24,[3^1,4^1,3^1,4^2]),\\
& (12,[3^1,4^4]),  (42,[3^1,7^1,3^1,7^1]), (24,[3^1,8^1,3^1,8^1]), (18,[3^1,9^1,3^1,9^1]), \\ 
& (84,[3^1,4^1,7^1,4^1]), (48,[3^1,4^1,8^1,4^1]), (36,[3^1,4^1,9^1,4^1]), (30,[3^1,5^3]), \\
&  (30,[3^1,4^1,10^1,4^1]),  (24,[3^1,6^1,4^1,6^1]), (15,[3^1,6^1,5^1,6^1]), (40,[4^3,5^1]), \\
&  (24,[4^3,6^1]), (20,[5^1,4^1,5^1,4^1]),(168,[4^1,6^1,14^1]), (96,[4^1,6^1,16^1]), \\
&  (72,[4^1,6^1,18^1]),(60,[4^1,6^1,20^1]),(80,[4^1,8^1,10^1]), (48,[4^1,8^1,12^1]),\\
&   (84,[6^2, 7^1]),(48,[6^2, 8^1]), (36,[6^2, 9^1]), (30,[6^2,10^1]),(40,[5^1, 8^2]),\\
&  (24,[6^1, 8^2]), (40,[4^1, 10^2]), (24,[4^1, 12^2]),(84,[3^1, 14^2]),\\
&  (48,[3^1, 16^2]), (36,[3^1, 18^2]), (28,[7^3])\}.
\end{align*}

Clearly, \begin{align*}
(n, X) = & (12, [3^7]), (24,[3^5,4^1]), (12,[3^4,4^2]), (12,[3^3,4^1,3^1, 4^1]),  (15,[3^5,5^1]),  (84,[3^4,7^1]), \\
&  (48,[3^4,8^1]), (36,[3^4,9^1]),  (30,[3^4,10^1]),  (20,[3^2,5^1,3^1,5^1]), (24,[3^3,4^1,6^1]),   \\
&  (24,[3^2,4^1,3^1,6^1]),  (24,[3^1,4^1,3^1,4^2]), (12,[3^1,4^4]), (42,[3^1,7^1,3^1,7^1]),    \\  
& (84,[3^1,4^1,7^1,4^1]), (48,[3^1,4^1,8^1,4^1]),  (40,[4^3,5^1]), (36,[3^1,4^1,9^1,4^1]),  \\
& (30,[3^1,4^1,10^1,4^1]), (30,[3^1,5^3]),  (24,[3^1,6^1,4^1,6^1]),  (72,[4^1,6^1,8^1]), (24,[4^3,6^1]),  \\
& (20,[5^1,4^1,5^1,4^1]), (168,[4^1,6^1,14^1]), (96,[4^1,6^1,16^1]), (60,[4^1,6^1,20^1]),   \\
& (80,[4^1,8^1,10^1]),  (48,[4^1,8^1,12^1]), (84,[6^2,7^1]),  (48,[8^1, 6^2]), (36,[9^1, 6^2]),  \\
& (30,[10^1, 6^2]), (40,[5^1, 8^2]),  (24,[6^1, 8^2]), (40,[4^1, 10^2]),  (84,[3^1, 14^2]) ~or~  (28,[7^3])
\end{align*} by Table 1 and Lemma \ref{lem2}. This completes the proof.
\end{proof}

\begin{proof}[Proof of Theorem \ref{theo:sf1}]
	Let $ \mathcal{K} $ be a map of type $ [3^4,10^1] $ on the surface with \Echar{-2}. Let $ V(\mathcal{K}) $, $ F^{10}(\mathcal{K}) $ denotes the set of vertices and set of 10-gon faces, respectively, of $\mathcal{K}$. Then, from the Euler characteristic equation, we have {$\# V(\mathcal{K})$}$ =30 $ and {$\# F^{10}(\mathcal{K})$}$ =3 $. With out loss of generality, assume $ V(\mathcal{K})=\{1,2,\dots, 30\} $ and $ F^{10}(\mathcal{K})=\{F_1,F_2,F_3\} $ with $ F_1=[1,2,\dots, 10] $, $ F_2=[21,22,\dots, 30] $, $ F_3=[11,12,\dots,20] $, and $ lk(1)=C_{12}(11,21,30,[2,3,4,5,\allowbreak6,7,8,9,10]) $. Then it is easy to see that the general graph for this type is as Figure \ref{3(4)_10_g}.
	\begin{figure}[H]
		\centering
		\includegraphics[height=7cm, width= 8cm]{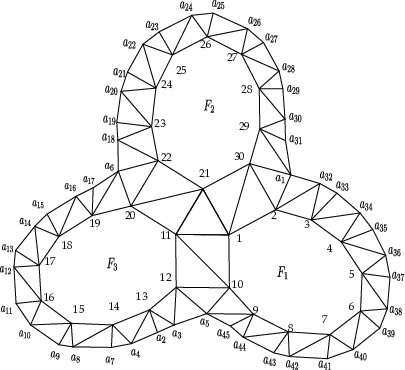}
		\caption{General graph of type $ [3^4,10^1] $}
		\label{3(4)_10_g}
	\end{figure}
	From Figure \ref{3(4)_10_g}, we can see that $[a_2,a_4]$, $[a_9,a_{10}]$, $[a_{13},a_{14}]$, $[a_6,a_{17}]$, $[a_{21},a_{22}]$, $[a_{25},a_{26}]$, $[a_{29},a_{30}]$ are adjacent with $ F_1$, $[a_3,a_5]$, $[a_7,a_8]$, $[a_{11},a_{12}]$, $[a_{15},a_{16}]$, $[a_{33},a_{34}]$, $[a_{37},a_{38}]$, $[a_{41},a_{42}]$, $[a_5,a_{45}]$ are adjacent with $ F_2 $, and $[a_{19},a_{20}]$, $[a_{23},a_{24}]$, $[a_{27},a_{28}]$, $[a_{1},a_{31}]$, $[a_1,a_{32}]$, $[a_{35},a_{36}]$, $[a_{29},a_{40}]$, $[a_{43},a_{44}] $ are adjacent with $ F_3$, and $ a_1\in \{14,15,16,17,18\} $.
	
	$ a_1=14 $ implies $ a_3\in \{23,24,25,26,27,28\} $. For $ a_3=23,25,27 $, these maps are same as $ \mathcal{K}_1[3^4,10^1] $, $ \mathcal{K}_2[3^4,10^1] $, and $ \mathcal{K}_3[3^4,10^1] $, respectively. For $ a_3=24,26,28 $, it is easy to see that $ [a_{21,a_{22}}]\in EG(F_3) $, $ [a_{25},a_{26}]\in EG(F_3) $ and $ [a_{29},a_{30}]\in EG(F_3) $, respectively, where $ EG(F_i) $ denotes the edge set of the face $ F_i $. This leads to a contradiction.
	
	For $ a_1=15,17 $, we get $ [a_2,a_4]\in EG(F_2) $ and $ [a_{13},a_{14}]\in EG(F_2) $, respectively. These are a contradiction.
	
	$ a_1=16 $ implies $ a_8\in \{23,24,\dots,28\} $. $ a_8=23 $ indicates that $ a_{16}\in \{25,26,27,28\} $. For $ a_{16}=25 $, we see that this map is same as $ \mathcal{K}_4[3^4,10^1] $. $ a_{16}=26,28 $ conclude that $ [a_{23},a_{24}]\in EG(F_1) $ and $ [a_{29},a_{30}]\in EG(F_3) $, respectively, which are contradictions. $ a_{16}=27 $ indicates $ a_{14}\in\{7,8,9\} $. $ a_{14}=7 $ conclude that $ [a_{41},a_{42}]\in EG(F_3) $, and $ a_{14}=8,9 $ implies that the face-sequence will not follow in $ lk(17) $, which are contradictions. $ a_8=24,26,28 $ indicates $ [a_{21},a_{22}] \in EG(F_3) $, $ [a_{19},a_{20}]\in EG(F_1) $ and $ [a_{29},a_{30}]\in EG(F_3) $, respectively. These make contradictions. $ a_8=25 $ implies $ a_4\in \{5,6,7,8\} $. For $ a_4=5,6 $, we see that $ [a_{37},a_{38}],[a_{29},a_{30}]\in EG(F_3) $, respectively. This is a contradiction. $ a_4=7 $  indicates either $ [a_{13},a_{14}]\in EG(F_2) $ or $ [a_{35},a_{36}]\in EG(F_2) $, which is a contradiction. $ a_4=8 $ implies $ a_3=a_{45} $, which is not possible.
	
	$ a_1=18 $ implies $ a_{12}\in \{23,24,\dots, 28\} $. $ a_{12}=24,26,28 $ indicates $ [a_{19},a_{20}]\in EG(F_1) $, $ [a_{25},a_{26}]\in EG(F_3) $ and $ [a_{29},a_{30}]\in EG(F_3) $, respectively, which is a contradiction. If $ a_{12}=23,27 $ then these maps are isomorphic with $ \mathcal{K}_3[3^4,10^1] $ and $ \mathcal{K}_4[3^4,10^1] $, respectively, under the map 
	\begin{align*}
	   &  (1, 11, 21)(2, 12, 22)(3, 13, 23)(4, 14, 24)(5, 15, 25) \\
	   &  (6, 16, 26)(7, 17, 27)(8, 18, 28)(9, 19, 29)(10, 2,0 30)
	\end{align*}
	
	 and if $ a_{12}=25 $ then this map is isomorphic with $ \mathcal{K}_1[3^4,10^1] $ under the map 
	 \begin{align*}
	     &  (1, 3, 5, 7, 9)(2, 4, 6, 8, 10)(11, 13, 15, 17 19)\\
	     &  (12, 14, 16, 18, 20)(21, 23, 25, 27, 29)(22, 24, 26, 28, 30) .
	 \end{align*}
	
	We call $ [a,b,c|d]$ a block in $ \mathcal{K} $ if $ [a,b,c] $ and $ [b,c,d] $ are two faces of $ \mathcal{K} $ and  the edges $ [a,b] $ and $ [c,d] $ are adjacent with 10-gonal faces. Observe that if a 10-gonal face intersect a block then the face adjacent with the block along an edge. Now, we define a cycle $C \colon \cdots \mbox{-} B_i \mbox{-}F_i\mbox{-}B_{i+1}\mbox{-}\cdots$ of blocks $B_i$ and 10-gonal faces $F_i$ such that $B_i, B_{i+1}$ are adjacent with $F_i$ antipodally. Again, we define a path $P_{B_1, B_m} \colon B_1\mbox{-}\cdots \mbox{-} B_i \mbox{-}\triangle_i\mbox{-} B_{i+1}\mbox{-}\cdots\mbox{-}B_m$ of blocks $B_i$ and triangular faces $\triangle_i$ such that $B_i, B_{i+1}$ are adjacent with $\triangle_i$. Distance between blocks (DBB) is the minimum number of blocks among in all the paths  between the blocks. 
	
	The list of all cycles of above type with DBB for each $ \mathcal{K}_i[3^4,10^1] $ for $ i=1,2,\dots,4 $, are listed bellow. (DBB presented for each consecutive blocks of the corresponding cycle)
	\begin{align*}
	  \mathcal{K}_1[3^4,10^1] \mbox{--} &   C_2 \colon F_3-[19,20,4|5]-F_1-[9,10,24|25]-F_2-[29,30,14|15],\\
	  & ~~~~~~DBB=\{2,2,2\}\\
	  \mathcal{K}_2[3^4,10^1] \mbox{--} & C_1 \colon F_3-[20,11,21|22]-F_2-[26,27,9|10]-F_1-[4,5,15|16],\\
	  &~~~~~~ DBB=\{1,3,3\}\\
	  & C_2 \colon F_3-[19,20,8|9]-F_1-[3,4,24|25]-F_2-[29,30,14,15], \\
	  & ~~~~~~DBB=\{2,1,2\}\\
	  & C_3 \colon F_3-[18,19,27|28]-F_2-[22,23,7|8]-F_1-[2,3,13|14],\\
	  & ~~~~~~DBB=\{1,3,2\}\\
	  & C_4 \colon F_3-[17,18,6|7]-F_1-[1,2,30|21]-F_2-[25,26,12|13],\\
	  & ~~~~~~DBB=\{3,1,3\}\\
	  & C_5 \colon F_3-[16,17,23|24]-F_2-[28,29,5|6]-F_1-[10,1,11|12],\\
	  & ~~~~~~DBB=\{4,2,2\}
	  \end{align*}
	  \begin{align*}
	  \mathcal{K}_3[3^4,10^1]  \mbox{--} & C_1 \colon F_3-[20,11,21|22]-F_2-[26,27,3|4]-F_1-[8,9,15|16],\\
	  & ~~~~~~ DBB=\{2,4,2\}\\
	  \mathcal{K}_4[3^4,10^1] \mbox{--} 
	  & C_1 \colon F_3-[20,11,21|22]-F_2-[26,27,7|8]-F_1-[2,3,15|16],\\
	  & ~~~~~~ DBB=\{3,4,1\}\\
	  & C_2 \colon F_3-[19,20,4|5]-F_1-[9,10,28|29]-F_2-[23,24,14|15],\\
	  & ~~~~~~ DBB=\{2,2,1\}\\
	  & C_3 \colon F_3-[18,19,25|26]-F_2-[30,21,1|2]-F_1-[6,7,13|14],\\
	  & ~~~~~~ DBB=\{2,2,4\}\\
	  & C_4 \colon F_3-[17,18,8|9]-F_1-[3,4,22|23]-F_2-[27,28,12|13],\\
	  & ~~~~~~ DBB=\{2,2,1\}\\
	  & C_5 \colon F_3-[16,17,29|30]-F_2-[24,25,5|6]-F_1-[10,1,11|12],\\
	  & ~~~~~~ DBB=\{3,3,1\}
	\end{align*}
	
	Clearly, by the above set of cycles and DBB sets, we can conclude that the maps $ \mathcal{K}_i[3^4,10^1] $ for $ i=1,2,3,4 $ are all non-isomorphic. This completes the proof.
	\end{proof}

\section{Semi-equivelar maps on higher genus surfaces}\label{covering}
 
There are three types of topological cycles in maps on surfaces. If a cycle does not bound a $2$-disk and not genus-separating, it is called {\em non-trivial}. If a cycle bounds a $2$-disk, then it is called {\em contractible}. If a cycle does not bound a $2$-disk but genus-separating, it is called {\em genus-separating}. 

\begin{lemma}\label{thm:SEM}
Let $\mathcal{K}$ be a semi-equivelar map of type $X$ on the surface {with} Euler char. $-2$. If $\mathcal{K}$ contains a non-trivial cycle then there exists a $2$-fold covering $\mathcal{K}^2$ of $\mathcal{K}$ of type $X$.
\end{lemma}

\begin{proof} Let $C$ be a non-trivial {shortest} (possibly!) cycle in $\mathcal{K}$ such that the cycle $C$ divides face-cycles of the vertices of $C$, that is, every sub-path $u\mbox{-}v\mbox{-}w \subset C$ of length two is a chord of the $face$-$cycle(v)$  at each vertex $v \in V(C)$. So, there are two sequences of faces $Y_1, Y_2, \dots, Y_k$ and $Z_1, Z_2, \dots, Z_{\ell}$ on two different sides of $C$ {and incident with $C$}. We cut $\mathcal{K}$ along the cycle $C$, and hence we get a map, namely, $\mathcal{K}_C$ which is bounded by two identical cycle $C$. We denote these boundary cycles by $C_Y$ and $C_Z$ where the faces $Y_1, Y_2, \dots, Y_k$ are incident with $C_Y$ and $Z_1, Z_2, \dots, Z_{\ell}$ are incident with $C_Z$ in $\mathcal{K}_C$. Let $C_Y := C(u_1, \dots, u_r)$ and $C_Z := C(w_1, \dots, w_r)$. Then, $C_Y$ identified with $C_Z$ by the map $u_i \rightarrow w_i$ for {all} $1 \le i \le r$ in $\mathcal{K}$, that is, $u_i = w_i$ for {all} $i$ in $\mathcal{K}$. So, $V(\mathcal{K}_C) = V(\mathcal{K} \setminus C) \cup V(C_Y) \cup V(C_Z)$ where $V(\mathcal{K} \setminus C) = \{v_1, v_2, \dots, v_m\}$. We consider another copy $\mathcal{K}_C'$ of $\mathcal{K}_C$ such that $\mathcal{K}_C \cong \mathcal{K}_C'$ where $V(\mathcal{K}_C') = \{v_1', v_2', \dots, v_m'\} \cup V(C_Y') \cup V(C_Z')$, $\partial \mathcal{K}_C' = C_Y' ~(= C(u_1', \dots, u_r')) \cup C_Z' ~(=C(w_1', \dots, w_r'))$,  $v \sim v'$, $u \sim u'$, $w \sim w'$ with the faces $Y_1', Y_2', \dots, Y_k'$ are incident with $C_Y'$ and $Z_1', Z_2', \dots, Z_{\ell}'$ are incident with $C_Z'$.

Since $f \colon \mathcal{K}_C \cong \mathcal{K}_C'$ by $v_j \rightarrow v_j'$, $u_i \rightarrow u_i'$ and $w_k \rightarrow w_k'$ for all $i, j, k$,  it follows that $Y_i \mapsto Y_i'$ and $Z_j \mapsto Z_j'$ by the map $f$.  We identify $C_Y$ with $C_Z'$ by the map $h_1 \colon u_i \rightarrow w_i'$ and $C_Z$ with $C_Y'$ by the map $h_2 \colon w_i \rightarrow u_i'$ for all $1 \le i \le r$. Hence, we get a map, namely, $\mathcal{K}^2(\mathcal{K}, \mathcal{K}') = \mathcal{K}_C \circhash^{h_1h_2} \mathcal{K}_C'$ of type $X$ with $\chi = -4$. This completes the proof.
\end{proof}

\begin{proof}[Proof of Theorem \ref{theo:sf2}] Let $\mathcal{K}$ be a semi-equivelar map of type $X$ on the surface of Euler char. $-2$. By Lemma \ref{thm:SEM}, there exists a semi-equivelar $\mathcal{K}^2$ of type $X$. We repeat same process for $m$ copies with same notations as follows. Let $f_{\ell, \ell+1} \colon \mathcal{K}_{C_{\ell}}^{\ell} \cong \mathcal{K}_{C_{\ell+1}}^{\ell+1}$ ($1 \le \ell \le m$ (addition in the suffix is modulo  $m $)) by $v_j^{\ell} \rightarrow v_j^{\ell+1}$, $u_i^{\ell} \rightarrow u_i^{\ell+1}$ and $w_k^{\ell} \rightarrow w_k^{\ell+1}$ for all $i, j, k$;  it follows that $Y_i^{\ell} \mapsto Y_i^{\ell+1}$ and $Z_j^{\ell} \mapsto Z_j^{\ell+1}$ by the map $f_{\ell, \ell+1}$.  We identify $C_Y^{\ell-1}$ with $C_Z^{\ell}$ by the map $h_{\ell-1, \ell} \colon u_i^{\ell-1} \rightarrow w_i^{\ell}$ and $C_Y^{\ell}$ with $C_Z^{\ell+1}$ by the map $h_{\ell, \ell+1} \colon u_i^{\ell} \rightarrow w_i^{\ell+1}$ for all $1 \le i \le r$. Hence, we get a map, namely, $\mathcal{K}^m(\mathcal{K}_1, \dots, \mathcal{K}_m) = \mathcal{K}_{C_1}^1\circhash\cdots\circhash \mathcal{K}_{C_m}^m$ of type $X$ with $\chi  = -2m$ for $m \ge 2$.
Thus, there is a semi-equivelar map of type $X$ from $\mathcal{K}$ with $\chi  = -2m$ for $m \ge 2$. Here, we present an application of the construction on one example. 

Let 
\begin{align*}
K:= &  \{[0,1,2],[0,1,11],[0,2,6],[0,6,8],[0,4,8],[0,4,10],[0,10,11],[1,2,3],\\
& [1,3,7],[1,7,9],[1,5,9],[1,5,11],[2,3,4],[2,4,8],[2,8,10],[2,6,10],\\
& [3,4,5],[3,5,9],[3,7,11],[3,9,11],[4,5,6],[4,6,10],[5,6,7],[5,7,11],\\
&[6,7,8],[7,8,9],[8,9,10],[9,10,11]\}
\end{align*}
(in \cite[(Section 2, $N_1$)]{DU2006}) which is a semi-equivelar map of type $[3^7]$ on the $2$-torus. The cycle $L = C_3(0, 6, 10)$ in $K$ is non-trivial. We cut $K$ along $L$. Hence, we get a map $Y$ with two boundary cycles $C_1, C_2$. We represent $(Y, C_{Y, 1}, C_{Y, 2})$ to be a map $Y$ with two boundary cycles $C_{Y, 1}, C_{Y, 2}$. Let $(K_i, C_{K_i, 1}, C_{K_i, 2})$ for $i=1, 2$ be two isomorphic copies of $(Y, C_{Y, 1}, C_{Y, 2})$, i.e., $K_i \cong Y, C_{K_i, 1} \cong C_{Y, 1}, C_{K_i, 2} \cong C_{Y, 2}$. Consider the map $Z := K_1 \circhash^{g_1g_2} K_2$ where $C_{K_1, 1}$ identified with $C_{K_{2}, 2}$ by $g_1 \colon C_{K_1, 1} \rightarrow C_{K_{2}, 2}$ $\&$  $C_{K_2, 1}$ identified with $C_{K_{1}, 2}$ by $g_2 \colon C_{K_2, 1} \rightarrow C_{K_{1}, 2}$ in $Z$.  Clearly, $Z$ is a semi-equivelar map of type $[3^7]$ of genus $3$. This $Z$ is $2$-fold cover of $K$. Again, consider three copies of $K$ and repeat the same construction as above with the same cycle. By repeating this process at each step, we get a semi-equivelar map of type $[3^7]$ with different genus. An example of a $3$-covering map is presented in Fig. \ref{fig:3-covering map} of $K$. 
\end{proof}

\section{Symmetric group of semi-equivelar maps}\label{symmetric}

\begin{figure}
    \centering
    \includegraphics[height=9cm, width= 10cm]{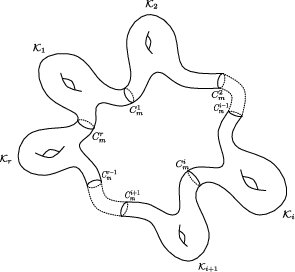}
    \caption{r-covering map}
     \label{fig:r-covering map}
\end{figure}

Let  $\mathcal{K}$ be a  map. Let $\mathcal{K}'$ be an isomorphic copy of $\mathcal{K}$ with $u \sim u'$ where $u \in V(\mathcal{K}), u' \in V(\mathcal{K}')$, that is, $f \colon \mathcal{K} \cong \mathcal{K}'$ where $f \colon V(\mathcal{K}) \mapsto V(\mathcal{K}')$ by $u \rightarrow u'$ (we used this notion throughout this section). Let $\mathcal{K}^2(\mathcal{K}, \mathcal{K}')$ be a $2$-covering map of $\mathcal{K}$ (as in Section \ref{covering}). In $\mathcal{K}^2(\mathcal{K}, \mathcal{K}')$, $\mathcal{K}, \mathcal{K}'$ identified along {non-trivial} cycles $C_m, C_m'$ (cycles of length $m$) where $C_m$ is in $\mathcal{K}$ and $C_m'$ is in $\mathcal{K}'$. Clearly, $C_m \sim C_m'$.
Similarly, let  $\mathcal{K}_1, \mathcal{K}_2, \dots, \mathcal{K}_r$ be $r$ isomorphic copies of $\mathcal{K}$.  Let $\mathcal{K}^r(\mathcal{K}_1, \mathcal{K}_2, \dots, \mathcal{K}_r)$ be a $r$-covering map of $\mathcal{K}$. Let $C_m^i$ be in $\mathcal{K}_i$ where $C_m^r \sim C_m$. {In  $\mathcal{K}^r(\mathcal{K}_1, \mathcal{K}_2, \dots, \mathcal{K}_r)$, $\mathcal{K}_1$, $\mathcal{K}_2$ identified along $C_m^1$; $\mathcal{K}_2$, $\mathcal{K}_3$ identified along $C_m^2$; \dots, and $\mathcal{K}_r$, $\mathcal{K}_1$ identified along $C_m^r$} (see in Fig. \ref{fig:r-covering map}).   

\begin{figure}
    \centering
    \includegraphics[height= 9cm, width= 10cm]{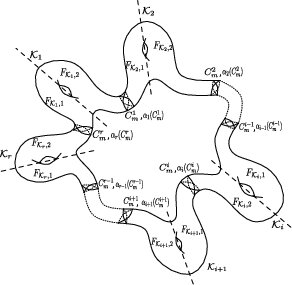}
    \caption{r-covering map-reflection}
     \label{fig:r-covering map-reflection}
\end{figure}

As it is discussed in Section \ref{covering}, there are two sequences of faces on two sides of $C_m$. We assume that the sequences of types of the faces are different. In fact, we can choose $C_m$ such that the assumption holds. Also assume the following: 
 let $\alpha \in Aut(\mathcal{K})$ be a product of disjoint transpositions such that $\alpha(F_{\mathcal{K}, 1}) = F_{\mathcal{K}, 2},$ $\alpha(C_m) \sim C_m$ where $F_\mathcal{K}$  = $F_{\mathcal{K}, 1} \cup F_{\mathcal{K}, 2}$ (call it ($*$)-property).
 Observe that if $C' (\neq C_m)$ is a cycle which is not homotopic to $C_m$ on the surface with $C' \sim C_m$  in $\mathcal{K}$ then $C_m \nsim C'$ in $\mathcal{K}^2$ (i.e., there does not exist any automorphism in $\mathcal{K}^2(\mathcal{K}, \mathcal{K}')$ which maps $C_m$ to $C'$) as $\mathcal{K}, \mathcal{K}'$ are identified in $\mathcal{K}^2$ along $C_m$. Similarly, if $\beta \in $Aut($\mathcal{K}$) where order of $\beta$ is greater than two ($o(\beta) >2$) and $\beta(C_m) \sim C_m$ then $\beta(C_m) \nsim C_m$ in  $\mathcal{K}^2$.
In general, let $F_{\mathcal{K}_t}  = F_{\mathcal{K}_t, 1} \cup F_{\mathcal{K}_t, 2}$ with $\alpha \sim \alpha_t, ~\alpha_t(F_{\mathcal{K}_t, 1}) = F_{\mathcal{K}_t, 2}$,  $1\le t \le r$. Then, for any other non-homotopic cycle $L$ in $\mathcal{K}^r(\mathcal{K}_1, \dots, \mathcal{K}_r)$, $C_m^i \nsim L (\neq C_m^i)$ for all $i$, and $\beta(C_m) \nsim C_m$ in $\mathcal{K}^r(\mathcal{K}_1, \dots, \mathcal{K}_r)$ if $o(\beta) >2$. See $C_m^i, \alpha_i(C_m^i), \mathcal{K}^r(\mathcal{K}_1, \dots, \mathcal{K}_r)$ in Fig. \ref{fig:r-covering map-reflection}. 

\begin{claim}\label{clm33}
Let $${\displaystyle \gamma= \prod_{ u \in V(F_{\mathcal{K}}), u_1 (\sim u) \in V(F_{\mathcal{K}_1}),\dots,~ u_r (\sim u) \in V(F_{\mathcal{K}_r})}(u_1, u_2, \dots, u_r)}.$$ Then, $\gamma \in $ Aut$(\mathcal{K}^r(\mathcal{K}_1, \mathcal{K}_2, \dots, \mathcal{K}_r))$.
\end{claim}

\begin{proof} Clearly, $F_{\mathcal{K}^r(\mathcal{K}_1, \mathcal{K}_2, \dots, \mathcal{K}_r)} = F_{\mathcal{K}_1} \cup F_{\mathcal{K}_2} \cup F_{\mathcal{K}_3} \cup \dots \cup F_{\mathcal{K}_r}$ and $\gamma(F_{\mathcal{K}_i})  = F_{\mathcal{K}_{i+1}}$ for all $i$  modulo $r$. Hence, $\gamma \in $ Aut$(\mathcal{K}^r(\mathcal{K}_1, \mathcal{K}_2, \dots, \mathcal{K}_r))$.
\end{proof}

\begin{claim}\label{clm34}
Let $${\displaystyle \gamma'= \prod_{u \in V(F_{\mathcal{K}}), u_1 (\sim u) \in V(F_{\mathcal{K}_1, 2} \cup F_{\mathcal{K}_2, 1}),\dots,~ u_r (\sim u) \in V(F_{\mathcal{K}_r, 2} \cup F_{\mathcal{K}_1, 1})}(u_1, u_2, \dots, u_r)}.$$ Then, $\gamma' \in $ Aut$(\mathcal{K}^r(\mathcal{K}_1, \mathcal{K}_2, \dots, \mathcal{K}_r))$ and $\gamma = \gamma'$.
\end{claim}

\begin{proof} In $\mathcal{K}_t$, $F_{\mathcal{K}_t}  = F_{\mathcal{K}_t, 1} \cup F_{\mathcal{K}_t, 2}$ such that $\alpha_t(F_{\mathcal{K}_t, 1}) = F_{\mathcal{K}_t, 2}$,  $1\le t \le r$. Clearly, $F_{\mathcal{K}^r(\mathcal{K}_1, \mathcal{K}_2, \dots, \mathcal{K}_r)} = (F_{\mathcal{K}_1, 2} \cup F_{\mathcal{K}_2, 1}) \cup (F_{\mathcal{K}_2, 2} \cup F_{\mathcal{K}_3, 1})\cup \dots \cup (F_{\mathcal{K}_r, 2} \cup F_{\mathcal{K}_1, 1})$ and $\gamma'(F_{\mathcal{K}_i, 2} \cup F_{\mathcal{K}_{i+1}, 1}) = F_{\mathcal{K}_{i+1}, 2} \cup F_{\mathcal{K}_{i+2}, 1}$ for all $i$  modulo $r$. Hence, $\gamma' \in $ Aut$(\mathcal{K}^r(\mathcal{K}_1, \mathcal{K}_2, \dots, \mathcal{K}_r))$. Observe that $\gamma'(F_{\mathcal{K}_i, 1}) =  F_{\mathcal{K}_{i+1}, 1}$ $\gamma'(F_{\mathcal{K}_i, 2}) =  F_{\mathcal{K}_{i+1}, 2},$ $\gamma(F_{\mathcal{K}_i, 1}) =  F_{\mathcal{K}_{i+1}, 1}$, $\gamma(F_{\mathcal{K}_i, 2}) =  F_{\mathcal{K}_{i+1}, 2}$ for all $i$  modulo $r$. Hence, $\gamma = \gamma'$.
\end{proof}

\begin{claim}\label{clm39}
  Aut$(\mathcal{K}^r(\mathcal{K}_1, \mathcal{K}_2, \dots, \mathcal{K}_r)) \cong \mathbb{D}_r$.
\end{claim}

\begin{proof}
 Let {$[a_1, a_2, \dots, a_r]$} be a $r$-gon. We know symmetry group of a polygon. The element of a symmetry group is either a rotation or a reflection symmetry.  Let $a_i = F_{\mathcal{K}_i, 2} \cup F_{\mathcal{K}_{i+1}, 1}$ or $F_{\mathcal{K}_i, 1} \cup F_{\mathcal{K}_i, 2}$  for all $i$  modulo $r$. So, the elements of Aut$(\mathcal{K}^r(\mathcal{K}_1, \dots, \mathcal{K}_r))$ are the product of rotations and reflections. 
 
  Let $S = \{F_{\mathcal{K}_1},\dots, F_{\mathcal{K}_r}\}$, that is, $a_i = F_{\mathcal{K}_i}$. Let $F_{\mathcal{K}_1},\dots, F_{\mathcal{K}_r}$ represent a {$[F_{\mathcal{K}_1}, \dots, F_{\mathcal{K}_r}]$}-gon. Let $X_j~ =~\prod_{A, B \in S}(A, B)~$ be a reflection symmetry of the polygon{,} $~I_{(A, B)}^j~ = \prod_{x \in A, y \in B, x \sim u, y \sim v, (u, v) \in \alpha} (x, y)${,} $Z_j = \{(a, b) \colon (a, b) \in I_{(A, B)}^j~where~(A,B) \in X_j \}${, and} $\delta_j = \prod_{(a,b) \in Z_j} (a, b)$. Observe that $\delta_j(F_{K^r}) = F_{K^r}$. Hence, $\delta_j \in $ Aut$(\mathcal{K}^r(\mathcal{K}_1, \dots, \mathcal{K}_r)).$ We know that the polygon has $r$ number of reflection symmetry, namely, $X_1, \dots, X_{r}$. Each $X_t$ defines $\delta_t$. Hence, $\mathcal{K}^r$ has the $\delta_1, \dots, \delta_{r}$ refection symmetry. Thus, we get the all possible reflection symmetry. Again, let  $a_i = F_{\mathcal{K}_i, 2} \cup F_{\mathcal{K}_{i+1}, 1}$.  We repeat {the} same process as above. Observe that in both the cases, we get the same list of reflection symmetry. In Fig. \ref{fig:3-covering map}, the map is an example of $3$-covering map of type $[3^7]$ of $N_1$ (in \cite[(Section 2, $N_1$)]{DU2006}); and 
  \begin{align*}
\delta_1 = & (0, 6)(10, 16)(2, 20)(1, 19)(7, 13)(8, 14)(9, 15)
(3, 21)(4, 22)(5, 23)\\
& (11, 17)(12, 30)(18, 24)(25, 31)(28, 34)(27, 33)(26, 32)(35, 29),\\
\delta_2 = & (12, 18)(22, 28)(14, 32)(13, 31)(19, 25)(20, 26)(21, 27)(15, 33)(16, 34)\\
& (17, 35)(23, 31)(24, 6)(30, 0)(1, 7)(4, 10)(3, 9)(2, 8)(11, 5),\\
\delta_3 = & (24, 30)(34, 4)(26, 8)(25, 7)(31, 1)(32, 2)(33, 3)(27, 9)(28, 10)(29, 11)\\
& (35, 7)(0, 18)(6, 12)(13, 19)(16, 22)(15, 21)(14, 20)(23, 17).
\end{align*}
  
  \begin{figure}
    \centering
    \includegraphics[height= 10cm, width= 14cm]{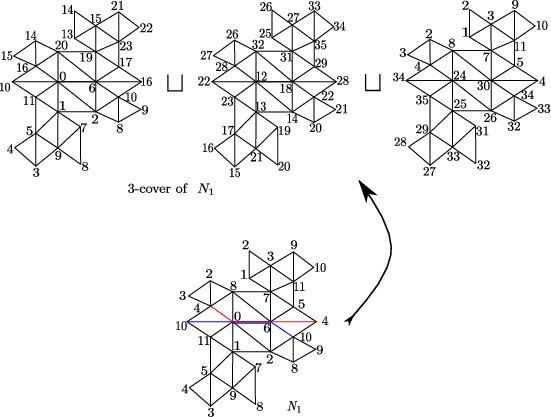}
    \caption{$3$-covering map of $N_1$ with type $[3^7]$ ($\alpha = (0,6)(4,10)(3, 9)(2, 8)(1, 7)(5, 11)$, $C_3(0, 6, 10)$, $\alpha(C_3(0, 6, 10)) = C_3(6, 0, 4)$)}
     \label{fig:3-covering map}
\end{figure}

 Let $a_i = F_{\mathcal{K}_i, 2} \cup F_{\mathcal{K}_{i+1}, 1}$. By Claim  \ref{clm33}, $\gamma$ is a rotation symmetry  in Aut$(\mathcal{K}^r(\mathcal{K}_1, \dots, \mathcal{K}_r))$. Similarly, if $a_i = F_{\mathcal{K}_i, 1} \cup F_{\mathcal{K}_{i}, 2}$, then, by Claim  \ref{clm34}, $\gamma'$ $\in$ Aut$(\mathcal{K}^r(\mathcal{K}_1, \dots, \mathcal{K}_r))$ and $\gamma = \gamma'$. Clearly, $C_m \nsim L$ for any cycle $L \neq C_m^1, \dots, C_m^r$ in $\mathcal{K}^r(\mathcal{K}_1, \dots, \mathcal{K}_r)$ as we identified $\mathcal{K}_1, \dots, \mathcal{K}_r$ along those cycles in $\mathcal{K}^r(\mathcal{K}_1, \dots, \mathcal{K}_r)$. Hence, there does not exist any other rotation symmetry in Aut$(\mathcal{K}^r(\mathcal{K}_1, \dots, \mathcal{K}_r))$.  In Fig. \ref{fig:3-covering map}, 
 \begin{align*}
\gamma = & (0, 12, 24)(1, 13, 25)(2, 14, 26)(3, 15, 27)(4, 16, 28)(5, 17, 29)(6, 18, 30)\\
& (7, 19, 31)(8, 20, 32)(9, 21, 33)(10, 22, 34)(11, 23, 35).
\end{align*}
 
 So, the symmetry group of $\mathcal{K}^r(\mathcal{K}_1, \mathcal{K}_2, \dots, \mathcal{K}_r)$ is generated by the rotation symmetry $\gamma$ ($o(\gamma) = r$) and the reflection symmetry $\delta_1, \delta_2, \dots, \delta_r$. Hence, Aut$(\mathcal{K}^r(\mathcal{K}_1, \mathcal{K}_2, \dots, \mathcal{K}_r))$ $= <\gamma, \delta_1, \delta_2, \dots, \delta_r> \cong \mathbb{D}_r$.   
\end{proof}


If $\alpha \not\in$ Aut($\mathcal{K}$) in above (($*$)-property). Then, we have the follow result.

\begin{claim}\label{clm40}
  Aut$(\mathcal{K}^r(\mathcal{K}_1, \mathcal{K}_2, \dots, \mathcal{K}_r)) \cong \mathbb{Z}_r$.
\end{claim}

\begin{proof} Assume that $\delta \in Aut(\mathcal{K}^r(\mathcal{K}_1, \dots, \mathcal{K}_r))$ is a reflection symmetry. Then $\delta(F_{\mathcal{K}_i}) = F_{\mathcal{K}_i}$ in $\mathcal{K}^r(\mathcal{K}_1, \dots, \mathcal{K}_r)$ for some $i$. This implies that $\alpha_i (\sim \alpha)$ exists in $\mathcal{K}_i$. This is a contraction as $\alpha \not\in$ Aut($\mathcal{K}$). So, there does not exist any reflection symmetry in Aut$(\mathcal{K}^r(\mathcal{K}_1, \mathcal{K}_2, \dots, \mathcal{K}_r))$. Thus, the symmetry group of $\mathcal{K}^r(\mathcal{K}_1, \dots, \mathcal{K}_r)$ is generated by the only rotation symmetry $\gamma$. It is clear from Claim \ref{clm39}. Hence, Aut$(\mathcal{K}^r(\mathcal{K}_1, \mathcal{K}_2, \dots, \mathcal{K}_r)) = <\gamma> \cong \mathbb{Z}_r$.  
\end{proof}

\begin{proof}[Proof of Theorem \ref{theo:sf3}] The proof follows from Claims \ref{clm39} and \ref{clm40}.
\end{proof}

{\small

}

\end{document}